\documentclass{amsart}

\usepackage{tabu}
\usepackage[margin=95pt]{geometry}
\usepackage{enumitem}
\allowdisplaybreaks

\usepackage{amssymb}
\usepackage{amsmath}
\usepackage{amscd}
\usepackage{amsbsy}
\usepackage{comment}
\usepackage{enumitem}
\usepackage[matrix,arrow]{xy}
\usepackage{hyperref}
\DeclareMathOperator{\Norm}{Norm}

\newcommand{\Q}{{\mathbb Q}}
\newcommand{\Z}{{\mathbb Z}}
\newcommand{\C}{{\mathbb C}}
\newcommand{\R}{{\mathbb R}}

\newcommand{\p}{\mathfrak{p}}

\def\mod#1{{\ifmmode\text{\rm\ (mod~$#1$)}
\else\discretionary{}{}{\hbox{ }}\rm(mod~$#1$)\fi}}

\begin {document}

\newtheorem{theorem}{Theorem}
\newtheorem{lemma}{Lemma}[section]
\newtheorem{prop}[lemma]{Proposition}

\theoremstyle{definition}

\theoremstyle{remark}

\title[Effective $S$-unit equations]
{Effective $S$-unit equations beyond $3$ terms : Newman's Conjecture}

\author[Prajeet Bajpai]{Prajeet Bajpai}
\address{Department of Mathematics, University of British Columbia, Vancouver, B.C., V6T 1Z2 Canada}
\email{prajeet@math.ubc.ca}

\author[Michael Bennett]{Michael A. Bennett}
\address{Department of Mathematics, University of British Columbia, Vancouver, B.C., V6T 1Z2 Canada}
\email{bennett@math.ubc.ca}

\thanks{The authors are supported by NSERC}

\date{\today}


\begin {abstract}
We show how to effectively solve $5$-term $S$ unit equations when the set of primes $S$ has cardinality at most $3,$ and use this to provide an explicit answer to an old question of D. J. Newman on representations of integers as sums of $S$-units.
\end {abstract}
\maketitle

\section{Introduction}\label{intro}

Let $K$ be a number field and $S$ a finite set of places of $K$, containing the infinite places $S_\infty$. We call $u$ an $S$-unit if $\lVert u \rVert_\nu = 1$ for all places $\nu$ of $K$ lying outside $S$. A quite remarkable collection of number theoretic problems, ranging from classical Diophantine equations, to the determination of representatives of binary forms of bounded discriminant, to monogenicity of number fields, to questions in transcendental number theory, may be reduced to that of solving  $S$-unit equations of the form
\begin{equation} \label{main-eq1}
a_1u_1+a_2u_2+\ldots+a_nu_n = 0, \; \; \; 
 \gcd (a_1 u_1, a_2 u_2, \ldots, a_n u_n)=1,
\end{equation}
where the $a_i$ are fixed nonzero elements of a number field $K$ and the $u_i$ are $S$-units of $K$. For excellent surveys of such applications, the reader is directed to the classic paper of Evertse, Gy\H{o}ry, Stewart and Tijdeman \cite{EGST} and to the more recent monograph of Evertse and Gy\H{o}ry \cite{EG}. While many finiteness results have been established for these equations, they are in general ineffective. In other words, for an arbitrary $S$-unit equation, it is often not possible, even in principle, to provide an upper bound on the height of potential solutions. 

Essentially all  effective results for $S$-unit equations are derived from Baker's bounds for linear forms in logarithms. Indeed, this method allows for an effective resolution of equation (\ref{main-eq1}) in case $n=3$ (see e.g.  \cite{BaDa} for an example where the $u_i$ are units, rather than $S$-units). In general, however, extending such results to the cases $n \geq 4$ remains, in almost every situation, beyond current technology.
Vojta, in his PhD thesis \cite{Voj} from 1983, however, was able to effectively solve equation (\ref{main-eq1}) if $n=4$ and $|S| \leq 3$,  via a pigeonhole argument (see also Skinner \cite{Sk} for an independent but less general result). While this is quite a strong restriction on the size of $S$ (in particular it only allows us to solve such equations in number fields of degree $\le 6$ and of particular signatures), this result nevertheless provided the first effective upper bounds for a number of Diophantine problems (see e.g. \cite{Sk},  \cite{TD}, \cite{TW} and \cite{Wang}).

In this paper, we extend this chain of reasoning to  equation (\ref{main-eq1}) with $n=5$, effectively solving equation (\ref{main-eq1}) in this case,
again under the restriction $|S| \leq 3$. 
The proof is once again based on a pigeonhole argument, albeit a somewhat subtle one that fully utilizes the  strength of bounds for linear form in logarithms. Let 
$\overline{u} = (u_1,u_2,u_3,u_4,u_5)$
be a solution to \eqref{eq-5!}. By applications of bounds for linear forms in complex and $p$-adic logarithms, we know that at least 3 of the $u_i$ must be {\it large} at any given place in $S$ -- meaning that they are of comparable size to the largest term at that place. Via the product formula (see Section \ref{prodform}), this implies that at least two terms in any solution of large height are in fact of comparable size at every place. In particular, if $u_1$ and $u_2$ are two such terms, we claim it is then possible to write $a_1u_1+a_2u_2 = a u$ with $u$ an $S$-unit, such that the new coefficient $a$ has small height compared to $\overline{u}$. We term this procedure {\it matching}. Equation (\ref{eq-5!}) then reduces to an $S$-unit equation with only four terms, which can be solved by arguing as in \cite{Voj}. Our main result is the following.
\begin{theorem}\label{eqmainthm}
Let $K$ be a number field and $S$ a finite set of places of $K$ containing all infinite places, and let $a_1,a_2,\ldots,a_5$ be fixed non-zero elements of $K$. If $S$ has at most three elements, then there is an effectively computable upper bound on the heights of nondegenerate solutions of the equation
\begin{equation} \label{eq-5!}
a_1u_1+a_2u_2+a_3u_3+a_4u_4+a_5u_5 = 0
\end{equation}
when the $u_j$ are taken to be $S$-units. Here, we call a solution to (\ref{eq-5!}) degenerate if $a_iu_i+a_ju_j=0$ for some $1 \leq i < j \leq 5$ (i.e. if equation (\ref{eq-5!})  has a {\it vanishing subsum}).
\end{theorem}

While this result may appear to be extremely restrictive, it and the argument behind it can actually be applied to effectively solve a variety of Diophantine problems. While analogues for simultaneous equations have been applied by the first author \cite{Baj23} and by the first author and Bugeaud \cite{BajBu}, we will concentrate here on applications of Theorem \ref{eqmainthm} and the machinery behind it. Specifically, we will apply our matching procedure to address a problem raised by Erd\H{o}s and Graham in \cite{ErGr}.  Let $N$ be a positive integer and write $\omega(N)$ for the number of distinct nonnegative integer tuples $(a,b,c,d)$ such that
\begin{equation}\label{newmanrep}
N = 2^a3^b + 2^c + 3^d .
\end{equation}
Here, we call two such representations $(a_1,b_1,c_1,d_1)$ and $(a_2,b_2,c_2,d_2)$ \emph{distinct} if 
\[
\{ 2^{a_1}3^{b_1}, 2^{c_1}, 3^{d_1}  \} \neq \{ 2^{a_2}3^{b_2}, 2^{c_2}, 3^{d_2}  \}.
\]
An old question of D.J. Newman (see Erd\H{o}s and Graham, page 80 of \cite{ErGr}) is whether $\omega(N)$ is absolutely bounded. This was settled in the affirmative by Evertse, Gy\H{o}ry, Stewart and Tijdeman in \cite{EGST}, with their argument subsequently refined by Tijdeman and Wang \cite{TW}, who proved
\begin{theorem}[Tijdeman and Wang, 1988]\label{TW} There exists a constant $N_0$ such that $\omega(N)\le 4$ for all $N>N_0$ (and hence a second constant $\omega_0$ such that $\omega(N)<\omega_0$ for all $N\ge 1$). 
\end{theorem}

The proof of this result is, however, ineffective, in that there is no way to extract an explicit value for $N_0$ or a sharp bound for $\omega_0$ in Theorem \ref{TW} from their arguments. We will derive versions of these results with explicit constants. 
Specifically, we will prove 
\begin{theorem} \label{thm-New}
Let $N$ be a positive integer. Then 
$$
\omega (N) \leq 
\begin{cases}
4 \; \; \mbox{ if }  \; \; N \ge 131082, \\
5 \; \; \mbox{ if }  \; \; N \ge 19700, \\
6 \; \; \mbox{ if }  \; \; N \ge 2316, \\
7 \; \; \mbox{ if }  \; \; N \ge 786, \\
8 \; \; \mbox{ if }  \; \; N \ge 300, \\
9 \; \; \mbox{ if }  \; \; N \ge 1. \\
\end{cases}
$$
We have $\omega(N)=9$ precisely when
\[
N \in \{ 41, 83, 89,113, 137, 161,227,299 \}.
\]
The largest $N$ with $\omega(N) = 5, 6, 7$ and $8$ are given by $N=131081, 19699, 2315$ and $785$, respectively. There are infinitely many $N$ with $\omega (N)=4$, corresponding to the identities
\begin{equation} \label{four-triv}
2^a+3^b = 2^{a-1} 3^0 + 2^{a-1}+3^b=2^{a-2} \cdot 3 + 2^{a-2} +3^b=2 \cdot 3^{b-1} + 2^a + 3^{b-1} = 2^3 3^{b-2} + 2^a + 3^{b-2}.
\end{equation}
\end{theorem}

The outline of this paper is the following. In Section \ref{secmainthm}, we present our primary technical tools from the theory of linear forms in logarithms, and use them to prove a slight generalization of  Theorem \ref{eqmainthm}. In Section \ref{simul}, we briefly discuss analogous results for simultaneous equations. Sections \ref{newman}, \ref{newman1}, \ref{newman2} and \ref{newman3} are devoted to the proof of Theorem \ref{thm-New}, containing bounds for linear forms in two logarithms, complex and $p$-adic, specialized to Newman's problem, results on $5$-term $S$-unit equations, a treatment of cases with vanishing subsums, and the final proof, respectively.

\section{Five-Term $S$-Unit Equations} \label{secmainthm}

\subsection{Technical Preliminaries}\label{prodform}
Let $K$ be a number field with $[K:\Q]=D$ and $\nu$ a place (equivalence class of absolute values) of $K$. If $\nu$ is archimedean, corresponding to an embedding $\sigma$ of $K$ in $\C$, we define
\[
\lVert x \rVert_\nu = \begin{cases}
|\sigma(x)| & \text{if } K_\nu = \R\\
|\sigma(x)|^2 & \text{if } K_\nu = \C,
\end{cases}
\]
while if $\nu$ is non-archimedean with associated rational prime $p$,  we normalize $\lVert\; \cdot \;\rVert_\nu$ so that
\[
\lVert p \rVert_\nu = p^{-[K_\nu:\Q_p]}.
\]
With these choices, we have the \emph{product formula}
\[
\prod_\nu \;\lVert x \rVert_\nu = 1,
\]
valid for all non-zero $x$ in $K$, where $\nu$ ranges over places corresponding to each prime ideal in $\mathcal{O}_K$ and to each real embedding and pair of complex conjugate embeddings (for details, see for example \cite{bombieri}). We let $S_\infty$ denote the set of archimedean places of $K$.

If $\overline{x} = (x_0,x_1,\ldots,x_n)\in \mathbb{P}^n(K)$, we define the Weil height
\[
H(\overline{x}) = \prod_{\nu}\max\{ \lVert x_0  \rVert_\nu, \lVert x_1  \rVert_\nu, \ldots , \lVert x_n  \rVert_\nu \}
\]
with associated logarithmic height $h(\overline{x}) = \log(H(\overline{x}))$. For $\alpha\in K$, we define the absolute logarithmic height of $\alpha$ via
\[
\frac1D\left( \log|a| + \sum_{i=1}^{D}\log\max\{1,|\alpha^{(i)}|\}  \right),
\]
where $a$ is the leading coefficient of the minimal polynomial of $\alpha$ over $\Z$, and the $\alpha^{(i)}$ are the conjugates of $\alpha$.

We will have need of bounds for linear forms in complex and $p$-adic logarithms;  results of Matveev \cite{Mat} and Yu \cite{Yu} are applicable in the broadest generality.

\begin{theorem}[Matveev, 2000]\label{matveev}
Let $\alpha_1,\ldots,\alpha_n$ be non-zero elements of $K$ and $b_1,\ldots, b_n$ be integers such that
\[
\Lambda = |b_1\log \alpha_1 + \cdots + b_n\log \alpha_n| \neq 0.
\]
If $K\subseteq\R$ put $\varkappa=1$, else put $\varkappa=2$. Let $B = \max\{|b_1|,\ldots,|b_n|\}$ and suppose that
\[
A_j \ge \max\{Dh(\alpha_j),|\log\alpha_j|,0.16)\},\quad 1\le j\le n.
\]
Then we have
\[
\log|\Lambda| \ge -\frac1\varkappa\left( en \right)^\varkappa 30^{n+3}n^{3.5}D^2\log(eD)\log(eB)A_1\cdots A_n.
\]

\end{theorem}

\begin{theorem}[Yu, 1998]\label{yu}
Let $\alpha_1,\ldots,\alpha_n$ be non-zero elements of $K$. Denote by $\p$ a prime ideal of the ring $\mathcal{O}_K$ of integers in $K$, lying above the rational prime $p$. Let $b_1,\ldots,b_n$ be rational integers such that 
\[
\Lambda = \alpha_1^{b_1}\cdots \alpha_n^{b_n} - 1 \neq 0
\]
and set $B\ge \max\{|b_1|,\ldots|b_n|,3\}$ and $A_j = \max\{h(\alpha_j),\log p\}$. Then we have
\[
\mathrm{ord}_\p(\Lambda) < 12\left( \frac{6(n+1)D}{\sqrt{\log p}} \right)^{2(n+1)} (p^D-1)\log(e^5nD)A_1\cdots A_n\log B.
\]
\end{theorem}

We will say that a solution $\overline{u} = (u_1,\ldots,u_n)$ to equation (\ref{main-eq1})  has a \emph{vanishing subsum} if 
\[
\sum_{i\in I}a_iu_i = 0\quad \text{for some non-empty } I \subsetneq \{1,2,\ldots,n\},
\]
and term such solutions {\it degenerate} (in that they correspond to solutions of $S$-unit equations in fewer-than-$(n-1)$ terms).

A final observation, important for our argument, is that Lemma 2.6 of \cite{Voj} can be slightly generalized to allow the $a_{i}$ to vary, provided they remain small in height relative to $\overline{u}$. More precisely, we have the following lemma.
\begin{lemma}\label{largelemma}
Let $\nu\in S$ and $\kappa_1,\kappa_2$ be absolute positive constants. Suppose that  
$$
\overline{u} = (u_1,u_2,u_3,u_4,u_5)
$$
is a nontrivial solution in $S$-units to the equation
\[
a_1u_1+a_2u_2+a_3u_3+a_4u_4+a_5u_5 = 0,
\]
satisfying $\log h(\overline u) >1$ and such that the (possibly variable) coefficients $a_i\; (1\le i\le 5)$ satisfy
\[
H(a_{i})\le \kappa_1 h(\overline{u})^{\kappa_2},
\]
while the $S$-units $u_j\; (1\le j \le 5)$ are ordered so that
\[
\lVert u_1 \rVert_\nu \ge \lVert u_2 \rVert_\nu \ge \lVert u_3 \rVert_\nu \ge \lVert u_4 \rVert_\nu \ge \lVert u_5 \rVert_\nu.
\]
Then there exist positive constants $c_1$ and $c_2$, depending only on $S,K$, $\kappa_1$ and $\kappa_2$, such that 
\[
\lVert u_3 \rVert_\nu \ge c_1 \lVert u_1 \rVert_\nu \, e^{-c_2\log^2 h(\overline{u}) }.
\]
\end{lemma}
\begin{proof}
The proof is identical to that of Lemma 2.6 in \cite{Voj} as the arguments using linear forms in (complex and $p$-adic) logarithms go through virtually unchanged. We reproduce it here in the interests of  keeping our exposition self-contained.

\textbf{Case 1:}
Suppose $\nu$ is archimedean and let $|\cdot|$ be the complex absolute value with respect to the embedding associated to $\nu$. Set
\[
A = \max\{|a_1|, \ldots, |a_5|, |a_1|^{-1}, \ldots, |a_5|^{-1}\}
\]
and note that our assumption on the heights of the $a_i$ implies that $A \ll h(\overline u)^c$ for a suitable constant $c$. We have
\[
3A|u_3| \ge | a_1u_1 + a_2u_2 | = | a_2u_2|\left|  -\frac{a_1u_1}{a_2u_2} -1  \right|.
\]
If $|a_1u_1| \le 0.9|a_2u_2|$, then we trivially have the inequalities 
$$
|u_3|\gg A^{-1}|u_1| \gg |u_1|e^{-c\log h(\overline u)}
$$
 and hence conclude as desired -- after changing constants, the same relation holds if we replace $|\cdot|$ by $\lVert\,\cdot\,\rVert_\nu$. If instead $|a_1u_1| > 0.9|a_2u_2|$, then 
\begin{equation}\label{somelogs}
A|u_3| \gg A^{-2}|u_1|\cdot \log\left(-\frac{a_1}{a_2}\cdot \frac{u_1}{u_2}\right)
\end{equation}
for a fixed determination of the complex logarithm. Let $w$ be the number of roots of unity in $K$, and let $\varepsilon_1,\ldots,\varepsilon_n$ generate the $S$-units of $K$ modulo roots of unity. Write $u_1/u_2$ as $\zeta \alpha_1^{b_1}\ldots \alpha_r^{b_n}$ for rational integers $b_1,\ldots, b_n$ and a root of unity $\zeta$. It remains to consider the linear form
\[
\Lambda = b_1\log\varepsilon_1+\cdots+b_n\log\varepsilon_r + b_{n+1}(2\pi i/w) + \log(-a_1/a_2)
\] 
where we can eliminate the discrepancies between $\log(\varepsilon_i^{b_i})$ and $b_i\log\varepsilon_i$ for each $1\le i \le n$ by suitably adjusting $b_{n+1}$. We now apply Theorem \ref{matveev} to conclude that
\[
|\Lambda| \ge e^{-C\log(eB)A_{n+2}}
\]
where  
$$
B = \max\{b_1,\ldots,b_{n+1}\}, \; \; \; A_{n+2} = \max\{Dh(-a_1/a_2),|\log(-a_1/a_2)|,0.16 \}
$$
and $C$ depends on $n$, the degree of $K$ and the heights of the generators of the group of $S$-units (one may consult, for example, \cite{Haj}  to see how these heights can be bounded). We have that $B\ll h(u_1/u_2)\ll h(\overline{u})$, and, by our assumption on the heights of the $a_i$,  $A_{n+2}\ll \log h(\overline{u})$. The various implied constants can be absorbed into the single constant $C$ above, and so, inserting this bound for $|\Lambda|$ into the relation \eqref{somelogs} derived above, 
\[
|u_3| \gg |u_1|\cdot e^{-C\log^2 h(\overline{u}) - D\log h(\overline{u})} \gg |u_1|\cdot e^{-C'\log^2 h(\overline{u})},
\]
which implies the same relation (although with new constants) if $|\cdot|$ is replaced by $\lVert\,\cdot\,\rVert_\nu$.

\textbf{Case 2:}
Suppose $\nu$ is non-archimedean, so $\lVert\, \cdot \,\rVert_\nu$ satisfies the triangle inequality. We have
\[
\lVert u_3 \rVert_{\nu} \gg \left( e^{-c\log h(\overline u)} \right) \lVert u_1 \rVert_{\nu}\left\lVert \zeta \alpha_1^{b_1}\ldots \alpha_r^{b_n}\left(  -\frac{a_1}{a_2} \right) - 1 \right\rVert_{\nu},
\]
similar to the previous case. Note that $\lVert u_1 \rVert_\nu = \lVert u_2 \rVert_\nu$. By Theorem \ref{yu}, we again have
\[
\left\lVert\alpha_1^{b_1}\ldots \alpha_r^{b_n}\left(  -\zeta \frac{a_1}{a_2} \right) -1\right\rVert_{\nu} \ge e^{-C\log B\cdot A_{n+2}}
\]
with $B = \max\{|b_1|,\ldots,|b_n|,3\}$ and $A_{n+2} = \max\{h(-\zeta a_1/a_2),\log p\}$. Here $C$ depends on $n$, the degree of $K$, and on the rational prime $p$ associated to the place $\nu$. If $B\le 3$, there are only finitely many solutions to consider. Otherwise, we note that $B\ll h(\overline{u})$ and $A_{n+2}\ll \log (h(\overline{u}))$ and thus once again obtain the desired result.
\end{proof}

\subsection{Proof of Theorem \ref{eqmainthm}}

We may assume, without loss of generality, that $|S|=3$.
Let $S=\{ \nu_1,\nu_2,\nu_3 \}$ and suppose that $\overline{u} = (u_1,u_2,u_3,u_4,u_5)$ is a non-trivial solution to \eqref{eq-5!}. For each $1\le i\le 3$, let $(u^{(i)}_1,u^{(i)}_2,u^{(i)}_3,u^{(i)}_4,u^{(i)}_5)$ be permutations of $(u_1,u_2,u_3,u_4,u_5)$ satisfying
\[
\left\lVert u^{(i)}_1 \right\rVert_{\nu_i} \ge \left\lVert u^{(i)}_2 \right\rVert_{\nu_i} \ge \left\lVert u^{(i)}_3 \right\rVert_{\nu_i} \ge \left\lVert u^{(i)}_4 \right\rVert_{\nu_i} \ge \left\lVert u^{(i)}_5 \right\rVert_{\nu_i}, \quad  1\le i\le 3
\]
i.e. such that the $i$-th permutation lists the terms $u_1,\dots,u_5$ in nonincreasing order of size at $\nu_i$.

By Lemma 2.6 of \cite{Voj} (i.e. by repeating the proof of Lemma \ref{largelemma} with constant bounds on the $a_j$), we have that
\[
\left\lVert u^{(i)}_3 \right\rVert_{\nu_i} \ge c_i \left\lVert u^{(i)}_1 \right\rVert_{\nu_i} h(\overline{u})^{-d_i},
\]
for positive constants $c_i, d_i$, for each $1\le i \le 3$. Thus at least three terms are {\it large} at every place. Moreover, being {\it large} at all three places contradicts  the product formula unless solutions are absolutely bounded (as in \cite{Voj}).

Thus, in the remaining cases, at least four terms must be large at two places each, whence at least two terms are large at the same two places. Without loss of generality, we may assume that $u_j$ and $u_k$ are large at $\nu_1$ and $\nu_2$, for some $1\le j<k\le 5$.  We thus have
\begin{align*}
\left\lVert u_j \right\rVert_{\nu_1} \ge c_1 \left\lVert u^{(1)}_1 \right\rVert_{\nu_1} h(\overline{u})^{-d_1},\; \left\Vert u_j \right\rVert_{\nu_2} \ge c_2 \left\lVert u^{(2)}_1 \right\rVert_{\nu_2} h(\overline{u})^{-d_2}\\
\left\lVert u_k \right\rVert_{\nu_1} \ge c_1 \left\lVert u^{(1)}_1 \right\rVert_{\nu_1} h(\overline{u})^{-d_1},\; \left\Vert u_k \right\rVert_{\nu_2} \ge c_2 \left\lVert u^{(2)}_1 \right\rVert_{\nu_2} h(\overline{u})^{-d_2}
\end{align*}
which in turn implies the inequalities
\begin{align}
\begin{split}\label{samesize} 
c\cdot h(\overline{u})^{-d} \le\;\, \lVert u_j/u_k \rVert_{\nu_1} \le c^{-1}\cdot h(\overline{u})^d \\
c\cdot h(\overline{u})^{-d} \le\;\, \lVert u_j/u_k \rVert_{\nu_2} \le c^{-1}\cdot h(\overline{u})^d,
\end{split}
\end{align}
for appropriately chosen positive constants $c$ and $d$ that depend only on $S$. Since the product formula yields
\[
\lVert u_j/u_k \rVert_{\nu_3} = \big(\lVert u_j/u_k \rVert_{\nu_1} \lVert u_j/u_k \rVert_{\nu_2}\big)^{-1},
\]
we also have
\begin{align}\label{samesize3}
c^2 \cdot h(\overline{u})^{-2d} \le\;\, \lVert u_j/u_k \rVert_{\nu_3} \le c^{-2}\cdot h(\overline{u})^{2d}
\end{align}
and thus $u_j$ and $u_k$ are of comparable size at each place in $S$.

By the main theorem of Hajdu \cite{Haj}, it is possible to choose fundamental $S$-units $\varepsilon_1,\varepsilon_2$ of bounded height and such that the entries of the matrix
\[A^{-1} = 
\begin{pmatrix}
\log \lVert \varepsilon_1 \rVert_{\nu_1} & \log \lVert \varepsilon_2 \rVert_{\nu_1} \\
\log \lVert \varepsilon_1 \rVert_{\nu_2} & \log \lVert \varepsilon_2 \rVert_{\nu_2}
\end{pmatrix}^{-1}
\]
satisfy
\begin{equation}\label{matrix}
\left|\left(A^{-1}\right)_{ij}\right| \le 2^{2-s}(s-1)!(s-2)!\frac{6d^4}{\log d},
\end{equation}
where $d$ is the degree of $K$ over $\Q$. Write
\[
u_j = \zeta_j \varepsilon_1^{\alpha_1}\varepsilon_2^{\alpha_2} \quad \text{and} \quad u_k = \zeta_k \varepsilon_1^{\beta_1} \varepsilon_2^{\beta_2},
\]
with $\zeta_j$ and $\zeta_k$ roots of unity in $K$ and $\alpha_\ell, \beta_\ell \in\Z$. Further, set 
\[
\gamma_1 = \min\{\alpha_1,\beta_1\} \; \; \mbox{ and } \; \;  \gamma_2 = \min\{\alpha_2,\beta_2\},
\]
and write
$$
a = a_j\zeta_j\varepsilon_1^{\alpha_1-\gamma_1} \varepsilon_2^{\alpha_2-\gamma_2} +  a_k\zeta_k\varepsilon_1^{\beta_1-\gamma_1}\varepsilon_2^{\beta_2-\gamma_2}
\; \; \mbox{ and } \; \; 
u = \varepsilon_1^{\gamma_1}\varepsilon_2^{\gamma_2}.
$$
Since
\[
\max\{\alpha_\ell - \gamma_\ell, \beta_\ell - \gamma_\ell\} = |\alpha_\ell - \beta_\ell|,
\]
we have
\[
\log \lVert a \rVert_{\nu} \le \log\big(\lVert a_j+a_k \rVert_{\nu}\big) + \sum_{i=1,2}\Big(|\alpha_\ell-\beta_\ell|\log\max\{ 1, \lVert \varepsilon_i \rVert_{\nu} \} \Big).
\]
Finally, writing $\bf{v}$ for the column vector $(\alpha_1-\beta_1,\alpha_2-\beta_2)$ and $\bf{w}$ for the vector 
$$
(\log \lVert u_j/u_k \rVert_{\nu_1}, \log \lVert u_j/u_k \rVert_{\nu_2}),
$$
we have $\bf{v}$ $= A^{-1}\bf{w}$. Thus, recalling \eqref{samesize} and \eqref{matrix}, we see that
\[
\log \lVert a \rVert_\nu \ll \log h(\overline{u})
\]
for all $\nu\in S$. Here the implicit constant depends only on the number field $K$, the primes in $S$, and the heights of the $a_i$. 

Thus, by writing $a_1u_1 + a_2u_2 = au$, we have reduced to a four-term $S$-unit equation with heights of the coefficients {\it small} in terms of $h(\overline{u})$. We may therefore move to considering equations
\[
a_1u_1 + a_2u_2 + a_3u_3 + a_4u_4 = 0
\]
which satisfy $H(a_j) \le \kappa_1 h(\overline{u})^{\kappa_2}$  ($1\le j \le 4$) for suitably chosen absolute constants $\kappa_1,\kappa_2$.

Again, for each $1\le i\le 3$, let $(u^{(i)}_1,u^{(i)}_2,u^{(i)}_3,u^{(i)}_4)$ be permutations of the $u_j$ that list them in decreasing order of size at $\nu_i$. Applying Lemma \ref{largelemma} gives
\[
\lVert u^{(i)}_3 \rVert_{\nu_i} \ge c_1 \lVert u^{(i)}_1 \rVert_{\nu_i} \, e^{-c_2\log^3 h(\overline{u}) },
\]
for each $1\le i \le 3$, where  $c_1$ and $c_2$ are effectively computable positive constants. Now, since there are only three primes in $S$,  at least one term $u_j$ does not occur as $u^{(i)}_4$ for any $1 \le i \le 3$.  This term is consequently large at all three places. However the product formula implies
\[
1 = \prod_{i=1}^3 \lVert u_j \rVert_{\nu_i} \ge c_1^3  \left( \prod_{i=1}^3\lVert u^{(i)}_1 \rVert_{\nu_i} \right) e^{-3c_2 \log^3 h(\overline{u})} =  c_1^3   H(\overline{u}) e^{-3c_2 \log^3 h(\overline{u})},
\]
and hence that $H(\overline{u})$ is absolutely bounded. This completes the proof of Theorem \ref{eqmainthm}.

We remark that we have actually proven the following somewhat stronger result (which will be of value to us later).
\begin{theorem}\label{eqmainthm1}
Let $K$ be a number field, $S$ a finite set of places of $K$ containing all infinite places, and   
$\overline{u} = (u_1,u_2,u_3,u_4,u_5)$
be a tuple of $S$-units with $\log h(\overline u) >1$.  Let $\kappa_1$ and $\kappa_2$ be positive constants and suppose that $a_1,a_2,\ldots,a_5$
are non-zero elements of $K$ with 
$$
H(a_{i})\le \kappa_1 h(\overline{u})^{\kappa_2},
$$
for $1 \leq i \leq5$. If $S$ has at most three elements and we have that
\[
a_1u_1+a_2u_2+a_3u_3+a_4u_4+a_5u_5 = 0,
\]
where $a_iu_i+a_ju_j\neq 0$ for all pairs $i, j$ with  $1 \leq i < j \leq 5$, then there is an effectively computable upper bound upon $H(\overline{u})$, depending only on $K, S, \kappa_1$ and $\kappa_2$.
\end{theorem}


\section{Simultaneous equations} \label{simul}

In Section \ref{secmainthm}, we were able to use lower bounds for linear forms in logarithms to deduce that the three largest terms at each place must be of {\it comparable} size. For systems of simultaneous equations, one can reduce the number of terms in each equation by Gaussian elimination and thereby win more large terms at each place. Under suitable conditions, this again allows us to match units. We describe the general argument below, which is based on similar ideas as the proof of Theorem \ref{eqmainthm}.

Suppose we have a system of $m$ equations in $n$ units $u_1,\ldots u_n$, given by
\[
a_{i1}u_1 + a_{i2}u_2 + \cdots + a_{in} u_n, \qquad 1 \le i \le m,
\]
and let $A = (a_{ij})$ be the corresponding matrix of coefficients. Suppose also that for every place $\nu_i \in S$ we have
\begin{equation}\label{klarge}
\left\lVert u_k^{(i)} \right\rVert_{\nu_i} \ge c\left\lVert u_1^{(i)} \right\rVert_{\nu_i} h(\overline u)^{-d},
\end{equation}
i.e. that $k$ terms are large at every place. We know that if an $S$-unit $u$ satisfies the above inequality for $|S|-1$ places, then, via the product formula, it is necessarily small at the remaining place in $S$. In particular, if two $S$-units $u_1$ and $u_2$ are large at the same $|S|-1$ places, then they are essentially the same size at every place in $S$.

Write $\alpha_i$ for the number of places in $S$ for which $u_i$ fails to be large (in the sense of \eqref{klarge}). From the product formula (and using Vojta's argument), we may assume that none of the units is large at every place of $S$. Our assumption \eqref{klarge} thus implies that
\[
(n-k)|S| \ge \sum_{i=1}^n \alpha_i \, ,
\]
where the $\alpha_i$ are positive integers. If $\alpha_i = 1$ for at most $|S|$ choices of subscript $i$, then it follows that
\[
(n-k)|S| \ge 2(n-|S|) + |S| = 2n - |S|.
\]
In other words, it follows  that if
\begin{equation}\label{matchineq}
(n-k+1)|S| < 2n
\end{equation}
then at least $|S|+1$ units are small at only one place each, and consequently two units $u_j$ and $u_\ell$ ($j \neq \ell$) are large in the sense of equation \eqref{klarge} at precisely the same $|S|-1$ places. Once again, by the same argument as at the end of Section \ref{secmainthm}, we can write
\[
a_ju_j + a_\ell u_\ell = a u
\]
with $u$ a unit and with $\log \lVert a \rVert_\nu \ll \log h(\overline u)$ for all $\nu \in S$. Thus we have succesfully reduced the number of units in the system from $n$ to $n-1$.

Let us see how to apply this matching argument to the system $A = (a_{ij})$. Fix a place $\nu \in S$ and without loss of generality assume $\lVert u_1 \rVert_\nu \ge \lVert u_2 \rVert_\nu \ge \ldots \ge \lVert u_n \rVert_\nu$. If the leftmost $m$ columns of $A$ are linearly independent, then after row-reducing we may assume $A$ has the form
\[
A = 
\begin{pmatrix}
a_{11}  &  &  & * & \cdots & * \\
 &   \ddots &  & \vdots & \ddots & \vdots \\
 &    & a_{mm} & * & \cdots & * \\
\end{pmatrix},
\]
i.e. the leftmost $m$ columns yield a diagonal matrix with each $a_{ii}$  non-zero. In particular, so long as we have
\[
a_{ii}u_i + a_{i(m+1)} \neq 0 \qquad \text{for some }1 \le i \le m,
\]
then (by an application of bounds for  linear forms in logarithms) at least $m+2$ units are large at $\nu$ in the sense of \eqref{klarge}. If \emph{any} choice of $m$ columns of $A$ have rank $m$, then the above argument works for all places $\nu \in S$.  Criterion \eqref{matchineq} with $k = m+2$, then yields the matching inequality
\begin{equation}\label{matchineq}
(n-m-1)|S| < 2n.
\end{equation}
Thus if all the above conditions, including the matching inequality, are satisfied then we can guarantee that at least two units can be matched with each other. Our next result describes a situation where matching once is sufficient to guarantee that the resulting system can be solved effectively.
\begin{theorem}\label{simulteq}
Let $A = (a_{ij})$ be an $(n-4)\times n$ matrix of elements of a number field $K$ such that any $n-4$ columns of $A$ are linearly independent. Let $S$ be a set of places of $K$ containing all infinite places and satisfying
\[
|S| < 2n/3.
\]
Then the solutions in $S$-units to the system of equations
\[
a_{i1}u_1 + a_{i2}u_2 + \cdots + a_{in} u_n \qquad 1 \le i \le n-4.
\]
can be `effectively determined'.
\end{theorem}
\begin{proof}
Fix a place $\nu$ and assume without loss of generality that $\lVert u_1 \rVert_\nu \ge \lVert u_2 \rVert_\nu \ge \ldots \ge \lVert u_n \rVert_\nu$ and that $A$ has the form
\[
A = 
\begin{pmatrix}
a_{11}  &  &  & * & \cdots & * \\
 &   \ddots &  & \vdots & \ddots & \vdots \\
 &    & a_{mm} & * & \cdots & * \\
\end{pmatrix} .
\]
If $a_{ii}u_i + a_{i(m+1)} \neq 0$ for some $1\le i \le m$, then as above at least $n-2$ units are large at $\nu$. 

For any other place $\nu' \in S$, we rearrange the $u_i$ to be in decreasing order of size; correspondingly we must permute the columns of $A$ to match. The leftmost $n-4$ columns of this new matrix are again of rank $n-4$ (by assumption), thus can again be diagonalized. This means we can repeat the argument above at any $\nu'$ as well, and thus at least $n-2$ units are large at each place of $|S|$. The inequality $|S| < (2/3)n$ is then precisely our matching inequality \eqref{matchineq}, implying two units $u_j$ and $u_\ell$ can be matched.

We combine columns $j$ and $\ell$ of $A$ accordingly, to obtain a new column of coefficients which now depend on $\overline u$, but whose heights are `logarithmically bounded' in terms of $h(\overline u)$. Note that any $n-3$ columns of $A'$ contain among them at least $n-4$ columns of $A$, and thus have rank at least $n-4$ (in particular $A'$ still has rank $n-4$). Furthermore, matching could not have created a column of zeroes, since any two columns of $A$ are linearly independent.

Finally, we note that
\[
\Big( (n-1)-(n-4) - 2 \Big)|S| = |S| < 2n/3 < n - 1
\]
since $n>3$. In summary, we have that  the new matrix $A'$ satisfies all constraints of Theorem 2.13 of \cite{Voj}, with the small modification that the entries are now allowed to vary up to small height. Vojta's proof via linear forms in logarithms proceeds exactly the same for $A'$, thus showing that the set of all solutions can be effectively determined.

All that is left is to consider the case where $a_{ii}u_i + a_{i(m+1)} = 0$ for all $i$ -- in this case the columns $n-3,n-2.n-1,n$ yield a system of $n-4$ equations in $4$ variables. If $n=5$, this is a single equation in four $S$-units, and the condition $|S|<(2/3)n$ implies that $|S| \leq 3$, whence this sub-system can also be solved effectively (following \cite{Voj}). If $n\ge 2$, then one can row-reduce again to leave at most 3 non-zero terms in each row and thus we are only left with 3-term $S$-unit equations, which can all be solved effectively.
\end{proof}

One can try to apply matching to other systems of equations, under weaker constraints than $|S| < 2n/3$. Indeed, as long as the inequality $(n-m-1)|S|<2n$ holds, and there are no vanishing subsums, then matching can help reduce the number of units in the resulting system. However, we do not in general have a guarantee that $m$ independent equations remain after matching -- this was a crucial ingredient in the proof of Theorem \ref{simulteq}. If, in a particular application, $m$ independent equations can somehow be guaranteed after matching, we note that the inequality \eqref{matchineq} gets further biased in our favour allowing us to match again. This process can be repeated iteratively.

To see an application, let $\xi$ be an algebraic number of degree $\ge 6$ such that $K = \Q(\xi)$ is Galois and totally complex, and consider the norm-form equation
\[
\Norm_{K/\Q} (x_0 + x_1\xi + x_2\xi^2 + x_3\xi^3) = m
\]
to be solved in rational integers $x_0,x_1,x_2,x_3$. This norm-form equation can be reduced to a system of unit equations satisfying the conditions of Theorem \ref{simulteq}, thus in particular Theorem \ref{simulteq} implies that the heights of solutions can be effectively bounded. Such equations were considered, and effectively resolved following this argument, in \cite{Baj23}. Moreover, setting $X = \max\{|x_0|,\ldots,|x_3|\}$, a lower bound of the form
\[
\Norm_{K/\Q} (x_0 + x_1\xi + x_2\xi^2 + x_3\xi^3) \ge X^\delta
\]
for positive effective $\delta$ is demonstrated in \cite{BajBu}. Note that not all systems considered in these two papers satisfy the criteria of Theorem \ref{simulteq}. Nevertheless, various special properties of the specific equations in question could be gainfully exploited to allow matching to work.

We end this section by noting that there do exist systems of $S$-unit equations satisfying the matching inequality \eqref{matchineq} that nevertheless remain resistant to our methods. For example, one can consider the norm-form equation
\[
\Norm(x + y\sqrt2 + z\sqrt3 + w\sqrt{-7}) = 1.
\]
Any solution to this equation arises from a unit $u = x + y\sqrt2 + z\sqrt3 + w\sqrt{-7}$ lying in the ring of integers of the field $K =\Q(\sqrt2,\sqrt3,\sqrt{-7})$. Such a unit necessarily satisfies a system of $4$ linear equations with coefficients in $K$. Further, our matching inequality is also satisfied, since in this case $S$ consists of only the infinite places of $K$ so $|S| = 4 = n/2$. Thus we are guaranteed to be able to match units, and in fact \emph{all} the units can be matched in pairs to reduce from 8 terms to 4 in one step. However, if we assume $w=0$ then this matching results in several vanishing subsums and the rank of the system of unit equations drops from $4$ to $1$ after matching. The remaining unit equation (in four units and with $|S| = 4$) corresponds  to the norm-form equation
\[
\Norm(x + y\sqrt2 + z\sqrt3) = 1.
\]
This equation appears widely in the literature (see \cite{CoZa}, \cite{Sch72}, \cite{Schbook}); its effective resolution appears out of reach of current technology.

\section{Newman's Conjecture : preliminaries} \label{newman}

The remainder of this paper will be devoted to proving Theorem \ref{thm-New}. To do this, we will have need for somewhat more specialized bounds for linear forms in logarithms than those afforded by Theorems \ref{matveev} and \ref{yu}. Specifically, we will appeal to the following results.
The first is  a state-of-the-art lower bound for linear forms in two complex logarithms due to Laurent (a special case of Corollary 2 of \cite{Lau}) :

\begin{prop} \label{arch2}
Let $\alpha_1$ and $\alpha_2$ be positive, real, multiplicatively independent algebraic numbers, with $\log \alpha_1$ and $\log \alpha_2$ positive, and 
$$
D = \left[ \mathbb{Q} (\alpha_1, \alpha_2) : \mathbb{Q} \right].
$$
Let $b_1$ and $b_2$ be positive integers. Define
$$
\Lambda = \left| b_2 \log \alpha_2 - b_1 \log \alpha_1 \right|.
$$
Further, let 
$$
b' = \frac{b_1}{D \log A_2} + \frac{b_2}{D \log A_1},
$$
where 
$$
\log A_i \geq \max \left\{ h (\alpha_i), \frac{\log \alpha_i}{D}, \frac{1}{D} \right\}.
$$
 Here, if $\alpha$ is algebraic of degree $d$ over $\mathbb{Q}$, we define the {\it absolute logarithmic height} of $\alpha$ via
$$
h (\alpha) = \frac{1}{d} \left( \log |a| + \sum_{i=1}^d \log \max \{ 1, |\alpha^{(i)}| \} \right),
$$
where $a$ is the leading coefficient of the minimal polynomial of $\alpha$ over $\mathbb{Z}$, and the $\alpha^{(i)}$ are the conjugates of $\alpha$. It follows that we have
$$
\log | \Lambda | \geq -17.9 \, D^4 \left( \max \{ \log b' + 0.38, 30/D, 1 \} \right)^2 \log A_1 \log A_2.
$$
\end{prop}

The second is a  slightly simplified version of Theorem 2 of Bugeaud-p-adic ??

\begin{prop} \label{BL}
  Let $p$ be a prime, $x_1/y_1$ and $x_2/y_2$ be nonzero rational numbers, define $g$ to be the smallest positive integer with the property that both $\nu_p ( (x_1/y_1)^g-1)$ and $\nu_p ( (x_2/y_2)^g-1)$ are positive and assume that there exists a real number $E$ such that
  $$
  \nu_p ( (x_1/y_1)^g-1) \geq E > 1/(p-1).
  $$
  If we define
  $$
  \Lambda = \left( \frac{x_1}{y_1} \right)^{b_1} - \left( \frac{x_2}{y_2} \right)^{b_2},
  $$
  where $b_1$ and $b_2$ are positive integers, then we have
  $$
  \nu_p (\Lambda) \leq \frac{36.1 \, g}{E^3 ( \log p )^4} \; \left( \max \{ \log b' + \log (E \log p) + 0.4, 6 E \log p, 5 \} \right)^2 \log A_1 \log A_2,
  $$
  if $p$ is odd or if $p=2$ and $\nu_2  (x_2/y_2-1) \geq 2$. If $p=2$ and $\nu_2  (x_2/y_2-1) \leq 1$, 
  $$
  \nu_2 (\Lambda) \leq  208 \; \left( \max \{ \log b' +  0.04, 10 \} \right)^2 \log A_1 \log A_2.
  $$
  Here
  $$
  b' = \frac{b_1}{\log A_2} + \frac{b_2}{\log A_1}
  $$
  and $A_1, A_2 > 1$ are real numbers such that
  $$
  \log A_i \geq \max \{ \log |x_i|, \log |y_i|, E \log p \}, \; i = 1, 2.
  $$
  \end{prop}

For our purposes, we will use these two results  to derive a number of bounds specifically tailored to Newman's Conjecture.

\begin{prop} \label{LFL2}
If $\alpha$ and $\beta$ are two positive integers, with $\beta \geq 2$, then
\begin{equation} \label{2-3-gap}
\left| 2^\alpha - 3^\beta \right| > \max \{ 2^\alpha, 3^\beta \}  \cdot \exp \left( -20 \log^2 \beta \right).
\end{equation}
\end{prop}

\begin{proof}

Suppose that (\ref{2-3-gap}) fails to be satisfied for integers $\alpha, \beta$ with $\beta \geq 2$. A short calculation ensures that $\beta \geq 100$. We thus have, applying Lemma? of Smart, 
\begin{equation} \label{JH}
 \Lambda < 1.1 \left| 2^\alpha 3^{-\beta} - 1 \right|  = \frac{ 1.1 \left| 2^\alpha - 3^\beta \right| }{\max \{ 2^\alpha, 3^\beta \}} < 1.1 \exp \left( -20 \log^2 \beta \right),
\end{equation}
where
$$
\Lambda = \left| \alpha \log 2 - \beta \log 3 \right|.
$$
We thus have
$$
\left| \frac{\log 3}{\log 2} - \frac{\alpha}{\beta} \right| < \frac{1.6}{\beta \exp \left( 20 \log^2 \beta \right)}.
$$
From this, necessarily $\alpha/\beta$ is a convergent in the infinite simple continued fraction expansion for $\log 3/\log 2$, say $\alpha/\beta = p_i/q_i$. Since we have
$$
\left| \frac{\log 3}{\log 2} - \frac{p_i}{q_i} \right| > \frac{1}{(a_{i+1}+2) q_i^2},
$$
where $a_{i+1}$ is the $(i+1)$-st partial quotient in the infinite simple continued fraction expansion for $\log 3/\log 2$, it follows that we necessarily have
$$
a_{i+1} > 0.6 \exp \left( 20 \log^2 \beta \right) \beta^{-1} -2 > 10^{181},
$$
since $\beta \geq 100$. A routine computation then implies that $i \geq 10000$ (the largest $a_{i+1}$ with $i < 10000$ is given by $a_{4313}=8228$), so that 
\begin{equation} \label{big}
\alpha > \beta \geq q_{10000} >10^{5144}.
\end{equation}

 We will apply Proposition \ref{arch2}  with
 $$
 b_2 = \alpha, \; \alpha_2 = 2, \; b_1=\beta, \; \alpha_1=3, \; D=1,
 $$
 so that $\log A_1= \log 3$, $\log A_2 = 1$, and
 $$
 b'= \beta + \frac{\alpha} {\log 3} < 2 \beta.
$$
 From (\ref{big}), we thus have
 $$
\log \Lambda \geq -17.9 \,  \left( 1.1+  \log \beta \right)^2 \log 3 > -19.8 \log^2 \beta,
$$
contradicting (\ref{JH}) and (\ref{big}).
  
\end{proof}

\begin{prop} \label{LFLP}
If $\gamma_1$ and $\gamma_2$ are odd integers, $u>10^5$ is an integer, and we set
$$
\gamma = \max \{ |\gamma_1|, |\gamma_2|, 3 \},
$$
then either $3^u \gamma_1 = \gamma_2$ or we have
$$
  \nu_2 (3^u \gamma_1 - \gamma_2 ) < 17 \log^2 u \log \gamma.
$$
\end{prop}
\begin{proof}
Suppose that $u > 10^5$, that
\begin{equation} \label{marsh}
  \nu_2 (3^u \gamma_1 - \gamma_2 ) \geq 17 \log^2 u \log \gamma
\end{equation}
and that $3^u \gamma_1 \neq \gamma_2$. In case $\gamma=3$, we have
$$
|3^u \gamma_1 - \gamma_2| \in \{ 3^u \pm 1, 3^{u+1} \pm 1, 3^{u} \pm 3, 3^{u+1} \pm 3  \}.
$$
As is well-known, if
$$
3^x \equiv 1 \mod{2^y},
$$
for $y \geq 2$, then $x \equiv 0 \mod{2^{y-2}}$ and hence, in particular, if $x > 0$, 
\begin{equation} \label{easy}
\nu_2(3^x-1) \leq 2 + \frac{\log x}{\log 2}.
\end{equation}
Further, the congruence $3^x \equiv -1 \mod{2^y}$ has no solutions for $y \geq 3$.
If $\gamma=3$, it follows that
$$
17 \log^2 (u) \log 3 \leq   \nu_2 (3^u \gamma_1 - \gamma_2 ) \leq 2 + \frac{\log (u+1)}{\log 2},
$$
contradicting $u \geq 2$. 

We may thus suppose that $\gamma \geq 5$. We make the observation that either  we have $\nu_2 \left( 3^u \gamma_1 - \gamma_2 \right)  =1$, or 
\begin{equation} \label{fungus}
\nu_2 \left( 3^u \gamma_1 - \gamma_2 \right)  = \nu_2 \left( 3^{2u}-\left(\gamma_2/\gamma_1 \right)^2 \right) -1.
\end{equation}
We may thus apply Proposition \ref{BL} with $p=2$, 
$$
x_1=9, \; y_1=1, \; x_2 = \pm \gamma_2, \; y_2= \gamma_1, \; g=1, \; E=3, \; b_1=u, \; b_2=2,
$$
where the sign of $x_2$ is chosen so that $x_2 \equiv y_2 \mod{4}$. We thus have
$$
\log A_1 = 2 \log 3 \; \mbox{ and } \; \log A_2 = \max \{ \log \gamma, 3 \log 2 \},
$$
whence
$$
  b' = \frac{u}{\max \{ \log \gamma, 3 \log 2 \}} + \frac{1}{\log 3}.
$$
If $\gamma \in \{ 5, 7 \}$, we thus have $\nu_2 (3^u \gamma_1 - \gamma_2 )$ bounded above by
$$
\frac{72.2 \log 3 }{9 ( \log 2 )^3} \; \left( \max \left\{ \log \left(  \frac{u}{3 \log 2} + \frac{1}{\log 3} \right) + \log (3 \log 2) + 0.4, 18 \log 2 \right\} \right)^2.
$$
Combining this with (\ref{marsh}) and using that $u > 10^5$, we find that
$$
17 \log^2 u \log \gamma < \frac{72.2 \log 3 }{9 ( \log 2 )^3} \; \left( \log \left(  \frac{u}{3 \log 2} + \frac{1}{\log 3} \right) + \log (3 \log 2) + 0.4 \right)^2.
$$
A short calculation reveals, in this latter case, that $u < 2.3 \times 10^{10}$ if $\gamma=5$, while $u \leq 42$ if $\gamma=7$, the latter being a contradiction.
  
  If we suppose that $\gamma=5$ then, from $u >10^5$ and (\ref{marsh}), necessarily
  $$
  \nu_2 (3^u \gamma_1 - \gamma_2) \geq 3200.
  $$
  This provides an immediate contradiction modulo $8$, if
  $$
  (\gamma_1,\gamma_2) \in \{ \pm (1,5), \pm (3,5), \pm (5,-5), \pm (5,1), \pm (5,3)  \}
  $$
  and contradicts (\ref{easy})  and the assumption that  $u < 2.3 \times 10^{10}$, if $ (\gamma_1,\gamma_2)  = \pm (5,5)$.
  It follows, then, that
   $$
  (\gamma_1,\gamma_2) \in \{ \pm (1,-5), \pm (3,-5),  \pm (5,-1), \pm (5,-3)  \}.
  $$
  To handle these remaining cases, note that
  $$
  3^x \equiv -5 \mod {2^{40}} \; \Rightarrow \;  x \equiv 205450132747 \mod{2^{38}}
  $$
  and
  $$
  5 \cdot 3^x \equiv -1 \mod {2^{40}} \; \Rightarrow \;  x \equiv 69427774197 \mod{2^{38}},
  $$
 again contradicting $u < 2.3 \times 10^{10}$. 

It follows, then, that $\gamma \geq 9$ and so we have
$$
  \nu_2 (3^u \gamma_1 - \gamma_2 ) < \frac{72.2 \log 3 }{27 ( \log 2 )^4} \; \left( \max \{ \log b' + \log (3 \log 2) + 0.4, 18 \log 2 \} \right)^2  \log \gamma,
$$
where now 
$$
  b' = \frac{u}{\log \gamma} + \frac{1}{\log 3}.
  $$
 If we have
$$
\max \{ \log b' + \log (3 \log 2) + 0.4, 18 \log 2 \} =18 \log 2,
$$
then
$$
  \nu_2 (3^u \gamma_1 - \gamma_2 ) < \frac{72.2 \times 12 \log 3 }{( \log 2 )^2} \;   \log \gamma < 17 \log^2 u \log \gamma,
$$
since $u > 10^5$. 
If, on the other hand, 
$$
\max \{ \log b' + \log (3 \log 2) + 0.4, 18 \log 2 \} = \log b' + \log (3 \log 2) + 0.4,
$$
we have
$$
17 \log^2 u \log \gamma < 13 \left( \log (u) + 0.3 \right)^2 \log \gamma,
$$
again contradicting $u > 10^5$.

\end{proof}

\begin{prop} \label{LFLP2}
If $\gamma_1$ and $\gamma_2$ are  integers with $\gcd (\gamma_1 \gamma_2,3)=1$,  $u>10^5$ is an integer, and we set
$$
\gamma = \max \{ |\gamma_1|, |\gamma_2|, 2 \},
$$
then either $2^u \gamma_1 = \gamma_2$ or we have
$$
  \nu_3 (2^u \gamma_1 - \gamma_2 ) < 10 \log^2 u \log \gamma.
$$
\end{prop}

\begin{proof}

Suppose that $u > 10^5$, 
\begin{equation} \label{marsh2}
  \nu_3 (2^u \gamma_1 - \gamma_2 ) \geq 10 \log^2 u \log \gamma
\end{equation}
and $2^u \gamma_1 \neq \gamma_2$. If $\gamma \leq 4$, then
$$
2^u \gamma_1 - \gamma_2 = 2^{u+\delta_1} \pm 2^{\delta_2}, \; \; \delta_i \in \{ 0, 1, 2 \}.
$$
If
$$
2^x \equiv 1 \mod{3^y},
$$
for $y \geq 1$, then $x \equiv 0 \mod{2 \cdot 3^{y-1}}$ and hence, if $x > 0$, 
\begin{equation} \label{easy2}
\nu_3(2^x-1) \leq 1 + \frac{\log (x/2)}{\log 3}.
\end{equation}
This immediately contradicts (\ref{marsh2}) if $\gamma \leq 4$. We may thus suppose that $\gamma \geq 5$.
We apply Proposition \ref{BL} with $p=3$,
$$
x_1=-8, \; y_1=1, \; x_2 = - \gamma_2, \; y_2= \gamma_1, \; g=1, \; E=2, \; b_1=u, \; b_2=3,
$$
We thus have
$$
\log A_1 = 2 \log 3 \; \mbox{ and } \; \log A_2 = \max \{ \log \gamma, 2 \log 3 \},
$$
whence
$$
  b' = \frac{u}{\max \{ \log \gamma, 2 \log 3 \}} + \frac{3}{2\log 3}.
$$
It follows that
$$
10 \log^2 u \log \gamma < \frac{36.1}{4  ( \log 3 )^3} \; \left( \max \{ \log b' + \log (2 \log 3) + 0.4, 12 \log 3 \} \right)^2  \log A_2.
$$
If $\gamma \in \{ 5, 7, 8 \}$, 
$$
10 \log^2 u \log \gamma < \frac{36.1}{2  ( \log 3 )^2} \; \left( \max \{ \log b' + \log (2 \log 3) + 0.4, 12 \log 3 \} \right)^2,
$$
so that either $u \leq 330433$ and  $\gamma=5$, or $u \leq 104543$ and $\gamma=7$, or we have 
$$
\max \{ \log b' + \log (2 \log 3) + 0.4, 12 \log 3 \}  = \log b' + \log (2 \log 3) + 0.4.
$$
In the last case, we have
$$
10 \log^2 u \log \gamma < \frac{36.1}{2  ( \log 3 )^2} \; \left( \log (u+3)+0.4 \right)^2,
$$
contradicting $\gamma \in \{ 5, 7, 8 \}$ and $u > 10^5$.
Arguing as in the proof of Proposition \ref{LFLP}, we contradict $10^5 < u \leq 330433$ (if $\gamma=5$) and $10^5 < u \leq 104543$ (if $\gamma=7$), by checking that the smallest positive integer $t$ for which one of 
$$
2^t \pm 5, \; 2^t \pm 7, \; 5 \cdot 2^t \pm 1,  \; 5 \cdot 2^t \pm 7, \; 7 \cdot 2^t \pm 1 \; \mbox{ or } 7 \cdot 2^t \pm 5
$$
is divisible by $3^{15}$ is $t=509305$.

If $\gamma \geq 10$, we have
$$
 10 \log^2 u  < \frac{36.1}{4  ( \log 3 )^3} \; \left( \max \{ \log b' + \log (2 \log 3) + 0.4, 12 \log 3 \} \right)^2,
$$
whence, again from $u > 10^5$, we have
$$
10 \log^2 u  < \frac{36.1}{4  ( \log 3 )^3} \; \left( \log (u+3)+0.4 \right)^2,
$$
once more contrary to $u > 10^5$.

\end{proof}

\section{Newman's Conjecture :  equations  with five or fewer terms} \label{newman1}

Fundamentally, Newman's Conjecture amounts to characterizing, for small integers  $k \geq 3$, solutions to the equation 
\begin{equation} \label{fund-k}
\sum_{i=1}^k u_i = 0,
\end{equation}
where the $u_i$ are rational numbers of the shape $\pm 2^\alpha 3^\beta$, for $\alpha$ and $\beta$ integers. Any such $k$-tuple corresponds to a unique $k$-tuple of integers $(u_1, \ldots, u_k)$ satisfying (\ref{fund-k}) with $P(u_i) \leq 3$,
\begin{equation} \label{order}
u_1 \geq u_2 \geq \cdots  \geq u_{k-1} \geq u_k, \; \; \gcd (u_1,u_2, \ldots, u_k)=1 \; \; \mbox{ and } \; \; u_1 > |u_k|.
\end{equation}
We will term such a solution to (\ref{fund-k})  {\it primitive}. As previously, if  a solution to (\ref{fund-k})  has the property that $\sum_{i \in I} u_i \neq 0$, for each proper nonempty subset $I$ of $\{ 1, 2, \ldots, k \}$, we will say that such a solution has {\it no vanishing subsums}.

A starting point for us are the following two results.

\begin{prop} \label{3-termite}
There are precisely $4$ primitive solutions to (\ref{fund-k}) in case $k=3$, with $P(u_i) \leq 3$, given by
$$
(u_1,u_2,u_3) = (2,-1,-1),  (3,-1,-2),  (4,-1,-3) \; \mbox{ and } \;  (9,-1,-8).
$$
\end{prop}

\begin{prop} \label{4-termite}
There are precisely $62$ primitive solutions to (\ref{fund-k}) in case $k=4$, with $P(u_i) \leq 3$ and no vanishing subsums, given by
$$
\begin{array}{l}
(u_1,u_2,u_3,u_4) = (3,-1,-1,-1), (3,1,-2,-2), (4,-1,-1,-2), (4,1,-2,-3), (4,2,-3,-3), \\
(6,-1,-2,-3),  (6,-1,-1,-4),  (6,1,-3,-4), (8,-2,-3,-3), (8,-1,-3,-4),  (8,-1,-1,-6), \\
 (8,1,-3,-6),  (9,-2,-3,-4), (9,-1,-4,-4),  (9,-1,-2,-6),  (9,1, -4,-6), (9,1,-2,-8),\\
 (9,2,-3,-8), (9,3,-4,-8), (12,-1,-3,-8), (12,-1,-2,-9),  (12,1,-4,-9), (16,-3,-4,-9),   \\
(16,-1,-6,-9),  (16,-1,-3,-12), (16,1,-8,-9),  (16,2,-9,-9), (18,-1,-8,-9),  \\
(18,-1,-1,-16),  (18,1,-3,-16), (24,1,-9,-16),  (27,-3,-8,-16),  (27,-2,-9,-16),  \\
(27,-1,-8,-18), 
 (27,-1,-2,-24),  (27,1,-12,-16), (27,1,-4,-24), (32,-2,-3,-27), \\
 (32, -1,-4,-27),   (32,1,-9,-24),  (32,1,-6,-27), (32,3,-8,-27),  (32,4,-9,-27), \\
  (36,-1,-8,-27), (36,-1,-3,-32),  (64,-1,-27,-36),  (64,-1,-9,-54),  (72,1,-9,-64), \\
 (81,-8,-9,-64),   (81,-1,-32,-48), (81,-1,-16,-64), (81,-1,-8,-72),  (81,1, -18,-64),   \\
  (96,1,-16,-81), (128,1,-48,-81), (144,1,-64,-81),   (256,-4,-9,-243), (256, -1,-12,-243), \\
 (256,3,-16,-243), (512,1,-27,-486), (512,1,-81,-432) \; \mbox{ and } \;  (729,-1,-216,-512). \\
\end{array}
$$
\end{prop}

The first of these is essentially a centuries-old result of Levi ben Gershon, while the second collects Theorems 1, 2 and 3, and Lemma 5 of Tijdeman and Wang \cite{TW}, Theorem 3 of Skinner \cite{Sk}, Theorems 1--5 of Wang \cite{Wang}, and Theorem 3.01 of Brenner and Foster \cite{BrFo}.

The main result of this section is the following.

\begin{theorem} \label{main-theorem}
There are precisely $1431$ primitive solutions to (\ref{fund-k}) in case $k=5$, with $P(u_i) \leq 3$ and no vanishing subsums.
The largest value of $u_1$ is
$$
u_1=3^{12}=531441,
$$
corresponding to one of the tuples
$$
\begin{array}{c}
(531441,2,-243,-6912,-524288), \; (531441,432,-1024,-6561,-524288) \\
\; \mbox{ or } \; (531441,-16,-576,-6561,-524288).
\end{array}
$$
\end{theorem}

The full list of these solutions is available at 

\hskip20ex \url{www.math.ubc.ca/~bennett/Newman-data}.

Equation (\ref{fund-k}), in case $k=5$, reduces to  $18$ families of equations of the shape
\begin{equation} \label{eq01}
2^a3^b+2^c+3^d+1=2^e3^f,
\end{equation}
\begin{equation} \label{eq02}
2^a3^b+2^e3^f+2^c+1=3^d,
\end{equation}
\begin{equation} \label{eq03}
2^a3^b+2^e3^f+3^d+1=2^c,
\end{equation}
\begin{equation} \label{eq04}
2^a+2^c+3^b+3^d=2^e3^f,
\end{equation}
\begin{equation} \label{eq05}
2^a+2^e3^f+3^b+3^d=2^c,
\end{equation}
\begin{equation} \label{eq06}
2^a+2^c+2^e3^f+3^b=3^d,
\end{equation}
\begin{equation} \label{eq07}
2^a3^b+2^c+3^d=2^e3^f+1,
\end{equation}
\begin{equation} \label{eq08}
2^a3^b+2^c+1=2^e3^f+3^d,
\end{equation}
\begin{equation} \label{eq09}
2^a3^b+3^d+1=2^e3^f+2^c,
\end{equation}
\begin{equation} \label{eq10}
2^c+3^d+1=2^e3^f+2^a3^b,
\end{equation}
\begin{equation} \label{eq11}
2^a3^b+2^e3^f+3^d=2^c+1,
\end{equation}
\begin{equation} \label{eq12}
2^a3^b+2^e3^f+2^c=3^d+1,
\end{equation}
\begin{equation} \label{eq13}
2^a3^b+2^e3^f+1=2^c+3^d,
\end{equation}
\begin{equation} \label{eq14}
2^a+2^c+3^b=2^e3^f+3^d,
\end{equation}
\begin{equation} \label{eq15}
2^a+3^b+3^d=2^e3^f+2^c,
\end{equation}
\begin{equation} \label{eq16}
2^a+3^b+2^e3^f=2^c+3^d,
\end{equation}
\begin{equation} \label{eq17}
3^b+3^d+2^e3^f=2^a+2^c
\end{equation}
and
\begin{equation} \label{eq18}
2^a+2^c+2^e3^f=3^b+3^d,
\end{equation}
where, in each case, the exponents are nonnegative integers.
For the final $12$ families (where there are two terms on the right-hand-side of the equation), to avoid vanishing subsums, we may suppose that no term on the left-hand-side equals one on the right-hand-side. 
Theorem \ref{main-theorem} is an almost immediate consequence of the following (noting that primitive solutions to (\ref{fund-k})  may correspond to several solutions to equations (\ref{eq01})--(\ref{eq18})).

\begin{theorem} \label{MT2}
Equations (\ref{eq01})--(\ref{eq18}) have precisely the following number of solutions without vanishing subsums, in tuples of nonnegative integers $(a,b,c,d,e,f)$ .
$$
\begin{array}{|cc|cc|cc|} \hline
\mbox{ Equation} & \# \mbox{ of solutions} & \mbox{ Equation} & \# \mbox{ of solutions}  & \mbox{ Equation} & \# \mbox{ of solutions}  \\ \hline
(\ref{eq01}) & 178 & (\ref{eq07}) & 151 & (\ref{eq13}) & 154 \\
(\ref{eq02}) & 127 & (\ref{eq08}) & 177 & (\ref{eq14}) & 298 \\
(\ref{eq03}) & 107 & (\ref{eq09}) & 149 & (\ref{eq15}) & 225 \\
(\ref{eq04}) & 298 & (\ref{eq10}) & 161 & (\ref{eq16}) & 219 \\
(\ref{eq05}) & 163 & (\ref{eq11}) & 119 & (\ref{eq17}) & 330 \\
(\ref{eq06}) & 191 & (\ref{eq12}) & 103 & (\ref{eq18}) & 302 \\ \hline
\end{array}
$$
In all cases, we have
$$
\max \{ a, c, e \} \leq 19 \; \; \mbox{ and } \; \; \max \{ b, d, f \} \leq 12.
$$
\end{theorem}

To prove this result, we treat each of these equations in turn. We will in fact provide full details for equation (\ref{eq07}) which will highlight all the main features of our arguments; analogous work for the other equations is available for the authors on request.

Suppose  that we have a solution to  equation (\ref{eq07}) in nonnegative integers $(a,b,c,d,e,f)$ with $cd \neq 0$ (which, in this case, is readily seen to be equivalent to the solution having no vanishing subsums). Define $M_1, M_2$ and $M_3$ such that
$$
\{ 3^d, 2^c, 2^a 3^b \} = \{ M_1, M_2, M_3 \},
$$
as multisets, with $M_1 \geq M_2 \geq M_3$. Let $M=2^e 3^f$ and suppose that
\begin{equation} \label{M-lower}
\log M > 3 \times 10^{10}.
\end{equation}
It follows that
$$
M_1 > \frac{1}{3} M.
$$
Applying  Proposition \ref{LFL2} to $M-M_1$,  
$$
M_2+M_3-1 > M \cdot \mbox{exp} \left( -20 \log^2 (\log M) \right),
$$
whence
\begin{equation} \label{eq-M2}
M_2 > \frac{1}{2} M   \cdot \exp \left( -20 \log^2 (\log M) \right).
\end{equation}
Since
$$
2^{\min \{ a, c, e \}} \mid 3^d-1
$$
and
$$
3^{\min \{ b, d, f \}} \mid 2^c-1,
$$
we deduce the inequalities
\begin{equation} \label{goodie1}
\min \{ a, c, e \}   \leq 2 + \frac{\log d}{\log 2}
\end{equation}
and
\begin{equation} \label{goodie2}
\min \{ b, d, f \}   \leq 1 + \frac{\log (c/2)}{\log 3}.
\end{equation}

Suppose first that  we have one of  $\min \{ b, d, f \} =d$, or $\min  \{ a, c, e \}  =c$, or both
\begin{equation} \label{fonzie!}
\min  \{ a, c, e \}  =a \; \mbox{ and } \min \{ b, d, f \} = b.
\end{equation}
Under these assumptions, it follows from (\ref{goodie1}) and (\ref{goodie2})  that
\begin{equation} \label{M3-upper}
M_3 < \frac{6 \log^2 M}{\log 2 \log 3}.
\end{equation}
Applying Propositions \ref{LFLP} and  \ref{LFLP2} to the forms
$$
M-(M_3-1), \; M_1+(M_3-1) \; \mbox{ and } M_2 + (M_3-1)
$$
 in each case contradicts (\ref{M-lower}), (\ref{eq-M2}) and (\ref{M3-upper}). By way of example, if we have (\ref{fonzie!}), then
$$
2^a 3^b \leq 6 c d < \frac{6 \log^2 M}{\log 2 \log 3},
$$
and hence, from (\ref{M-lower}) and (\ref{eq-M2}), $M_3=2^a 3^b$. Applying Proposition \ref{LFLP} to 
$$
3^d+(2^a3^b-1)
$$
 and Proposition \ref{LFLP2} to 
 $$
 2^c+(2^a3^b-1), 
 $$
 we find that 
  \begin{equation} \label{marcos7}
 \min \{ c, e \} < 17 \log^2 \left( \frac{\log M}{\log 3} \right)  \log \left(  \frac{6 \log^2 M}{\log 2 \log 3} \right)
 \end{equation}
 and
  \begin{equation} \label{marcos8}
 \min \{ d, f \} < 10 \log^2 \left( \frac{\log M}{\log 2} \right)  \log \left(  \frac{6 \log^2 M}{\log 2 \log 3} \right),
 \end{equation}
 contradicting $M=2^e 3^f$, $\min \{ 2^c,  3^d \} = M_2$, (\ref{M-lower}) and  (\ref{eq-M2}). We deduce similar contradictions if either $\min \{ b, d, f \} =d$ or $\min  \{ a, c, e \}  =c$ (and, in each case, require a rather weaker assumption than (\ref{M-lower})).

We may thus suppose that
$$
\left( \min \{ a, c, e \}, \min \{ b, d, f \} \right) = (e,b) \mbox{ or } (a,f).
$$
In these cases, we can rewrite one of $M-2^c$, $M-3^d$ or $M_1+M_2$ as $M M_0$, where $M_0$ is a positive rational number of {\it small} naive height. 
 To be precise, if say $\min \{ a, c, e \}=e$ and $M_3=2^a3^b$, we may apply Proposition \ref{LFLP}  to either
$$
(2^e3^{f-d}-1) 3^d +1  \; \mbox{ or } \; (2^e-3^{d-f}) 3^f +1
$$
to conclude from (\ref{eq-M2}) that
 \begin{equation} \label{marcos9}
 \min \{ a, c \} < 17 \log^2 \left( \frac{\log M}{\log 3} \right)  \left( \log 2 + \log^2 (\log M) \right). 
 \end{equation}
 Combining $\min \{ 2^c, 3^d \} = M_2$ with (\ref{M-lower}) and  (\ref{eq-M2}) implies that $ \min \{ a, c \} =a$ and hence, applying Proposition \ref{LFLP2}  to
$$
2^c + (2^a 3^b-1),
$$
\begin{equation} \label{marcos10}
 \min \{ d, f \} < 10 \log^2 \left( \frac{\log M}{\log 2} \right)   \left(  a \log 2 + b \log 3 \right).
 \end{equation}
With (\ref{M-lower}), (\ref{eq-M2}), (\ref{goodie2}) and (\ref{marcos10}), this is again a contradiction. We argue similarly if $\min \{ a, c, e \}=a$ or if $M_3 \neq 2^a3^b$.

It remains then to handle the values of $M$ with $\log M \leq 3 \times 10^{10}$, where we necessarily have
\begin{equation} \label{big}
\max \{ a, c, e \}  < 4.33 \times 10^{10} \; \mbox{ and } \; \max \{ b, d, f \} <  2.74 \times 10^{10}.
\end{equation}
This requires a reasonably short if somewhat unedifying case-by-case analysis. By way of example, if, say, 
$\left( \min \{ a, c, e \}, \min \{ b, d, f \} \right) = (e, b)$ and  $M_1=3^d$, it follows from (\ref{goodie1}) and (\ref{goodie2}) that $e \leq 36$ and $b \leq 22$. Since $3^d > 2^e 3^{f-1}$, we have either $d \leq f$, so that $d=f$ and $e=1$, or we have $1 \leq d-f \leq 22$. In the former case,  
$$
1+3^f = 2^c + 2^a 3^b,
$$
so that, from Proposition \ref{4-termite}, $(a,b,c,d,e,f)$ is one of
$$
\begin{array}{c}
(1,0,1,1,1,1), (1,0,3,2,1,2), (1,1,2,2,1,2), (1,2,6,4,1,4),  \\
(2,1,4,3,1,3), (3,0,1,2,1,2)  \mbox{ or }  (3,1,2,3,1,3).
\end{array}
$$
In the latter case, $2^e > 3^{d-f}$ and we have
\begin{equation} \label{final}
(2^e-3^{d-f}) 3^f = 2^c + 2^a 3^b-1,
\end{equation}
so that
$$
(2^e-3^{d-f}) 3^f \equiv -1 \mod{2^a}.
$$
For each pair $(e,d-f)$ with $2^e > 3^{d-f}$, $2 \leq e \leq 36$ and $1 \leq d-f \leq 22$, we check via Maple's ``ModularLog'' function  that the smallest nonnegative solution to the congruence
$$
(2^e-3^{d-f}) 3^f \equiv -1 \mod{2^{48}}
$$
corresponds to $e=26$, $d-f=1$ and $f=64541949951 >  2.74 \times 10^{10}$.  From (\ref{big}), it follows that $e \leq a \leq 47$. Looping over $e, d-f, a$ and $b$, we verify that the smallest positive $f$ for which
$$
(2^e-3^{d-f}) 3^f \equiv 2^a 3^b-1 \mod{2^{55}}
$$
corresponds to $e=7$, $a=35$, $d-f=4$ and $b=2$, where $f=247148451996 >  2.74 \times 10^{10}$.
We thus have $e \leq c \leq 54$ and, since $2^c \equiv 1 \mod{3^b}$, $b \leq 4$.

We now loop over $0 \leq b \leq 4$, $0 \leq a \leq 47$ and $1 \leq c \leq 54$ for which $2^c \equiv 1 \mod{3^b}$ and factor the right-hand-side of equation (\ref{final}) to find possible values for $f$. We then check to see if 
$$
\frac{2^c + 2^a 3^b-1}{3^f} + 3^{d-f}
$$
is, for some $1 \leq d-f \leq 22$, a power of $2$. In this special case, we find precisely the further solutions with $(a,b,c,d,e,f)$ one of 
$$
\begin{array}{c}
(2,1,4,4,2,3), (3,0,3,2,3,1), (3,1,2,4,2,3), (3,2,6,4,3,3), (4,1,4,4,4,2), \\
(5,0,5,4,4,2), (6,1,10,6,3,5), (6,1,10,8,5,5), (7,1,14,8,5,6)   \mbox{ or } (7,2,6,8,5,5).
\end{array}
$$

We treat the remaining cases in a similar fashion to arrive at Theorem \ref{MT2}; full details of our computations are again available  upon request.

\section{Newman's Conjecture : vanishing subsums } \label{newman2}

For the remainder of the paper, we will always suppose that a positive integer $N$ has the property that $\omega (N) \geq 2$; i.e. 
$$
N=2^{a_1}3^{b_1} + 2^{c_1} + 3^{d_1} = 2^{a_2}3^{b_2} + 2^{c_2} + 3^{d_2},
$$
where the $a_i, b_i, c_i$ and $d_i$ are nonnegative integers with 
\[
\{ 2^{a_1}3^{b_1}, 2^{c_1}, 3^{d_1}  \} \neq \{ 2^{a_2}3^{b_2}, 2^{c_2}, 3^{d_2}  \}.
\]
It follows that
\begin{equation} \label{varnish}
2^{a_1}3^{b_1} + 2^{c_1} + 3^{d_1} - 2^{a_2}3^{b_2} - 2^{c_2} - 3^{d_2}=0.
\end{equation}
We will say that such an equation has a {\it vanishing subsum} if 
$$
2^{a_1}3^{b_1} \in \{ 2^{a_2}3^{b_2}, 2^{c_2}, 3^{d_2}, 2^{a_2}3^{b_2} + 2^{c_2}, 2^{a_2}3^{b_2} + 3^{d_2}, 2^{c_2} +3^{d_2} \}.
$$
We will begin our proof of Newman's conjecture by defining a number of families of integers $N$ whose corresponding $S$-unit equations have {\it vanishing subsums}. In these cases, we will always be able to reduce equation (\ref{varnish}) to situations where we can appeal to the results of the previous section.

We will call an integer $N$  {\it type I special} if there exist nonnegative integers $a$ and $b$ with $N=2^a+3^b$. For these integers, the  identities 
(\ref{four-triv}) guarantee that $\omega (N) \geq 4$, at least provide $a, b \geq 2$. We say that an integer $N$ is {\it type II special} if there exist nonnegative integers $a$ and $b$ such that either
\begin{equation} \label{type2-1}
N=2^a+3^bc, \; \; c \in \{ 11, 19 \}, 
\end{equation}
or
\begin{equation} \label{type2-2}
N=2^a c+3^b, \;  \; c \in \{ 5, 7 \}.
\end{equation}
In these cases, the identities
$$
2^a+3^b \cdot 11 = 2 \cdot 3^{b} + 2^a + 3^{b+2}  = 2^3 \cdot 3^b + 2^a + 3^{b+1} = 2^5 \cdot 3^{b-1} + 2^a+ 3^{b-1},
$$
$$
2^a+3^b \cdot 19 = 2^4 \cdot 3^{b} + 2^a + 3^{b+1}  = 2 \cdot 3^{b+2} + 2^a + 3^{b} = 2^9 \cdot 3^{b-3} + 2^a+ 3^{b-3},
$$
$$
2^a \cdot 5 + 3^b = 2^{a+2} \cdot 3^0 + 2^a + 3^b = 2^{a} \cdot 3 + 2^{a+1} + 3^b = 2^{a-1} 3^2 + 2^{a-1} + 3^b
$$
and
$$
2^a \cdot 7 + 3^b = 2^{a+1} \cdot 3^1 + 2^a + 3^b = 2^{a} \cdot 3 + 2^{a+2} + 3^b = 2^{a-2} 3^3 + 2^{a-2} + 3^b
$$
ensure that $\omega (N) \geq 3$ for these $N$, if $a$ and $b$ are suitably large.

Finally, call an integer $N$ {\it type III special} if there exist nonnegative integers $a$ and $b$ with
\begin{equation} \label{type3-1}
N=2^a+3^bc, \; \; c \in \{ 5, 7, 13, 17, 25, 35, 43, 73, 97, 145, 259 \}, 
\end{equation}
\begin{equation} \label{type3-2}
N=2^a c+3^b, \;  \; c \in \{ 11, 13, 17, 19, 25, 35, 41, 73, 97, 145, 259 \},
\end{equation}
\begin{equation} \label{type3-3}
N=2^a 3^b+c, \;  \; b \in \{  1, 2 \} , \;  \; c \in \{ 3, 9 \},
\end{equation}
\begin{equation} \label{type3-4}
N=2^a 3^b+c, \;  \; a \in \{  1, 2 \} , \;  \; c \in \{ 2, 4 \},
\end{equation}
or
\begin{equation} \label{type3-5}
N=2^a 3^b+c,  \;  \; c \in \{ 5, 11, 17, 35, 259 \}.
\end{equation}

The main result of this section is the following.
\begin{theorem} \label{Vanishing}
If $N$ is a positive integer with $\omega (N) \geq 2$ and there exist two distinct representations of $N$ such that the corresponding equation (\ref{varnish}) has a vanishing subsum, then $N$ is special of type I, II or III. 

If $N$ is type I special, we have $\omega (N)\leq9$ for all $N$, $\omega (N)\leq4$ for all $N \ge 131082$ and $\omega (N)=4$, provided $N=2^a+3^b$ with $\min \{ a, b \} \geq 2$. The largest type I special $N$ with $\omega(N)=5, 6, 7, 8$ and $9$ are $131081, 19699, 2315, 283$ and $137$, respectively. 

If $N$ is type II special, then $\omega (N)\leq9$ for all $N$, $\omega (N) \le 3$ for all $N \ge 532308$ and $\omega (N)=3$ for all $N \ge 532308$ satisfying the additional constraints
$$
a\ge 1 \mbox{ if } c=5, \; \; a\ge 2 \mbox{ if } c=7, \; \;  b\ge 1 \mbox{ if } c=11 \; \mbox{ or } \; b\ge 3 \mbox{ if } c=19.
$$
The largest type II special $N$ with $\omega(N)=4, 5, 6, 7, 8$ and $9$ are $532307, 20483, 6665, 2267, 785$ and $299$, respectively.

If $N$ is type III special, then $\omega (N)\leq9$ for all $N$ and $\omega (N) \le 2$ for all $N \ge 76546076$. The largest type III special $N$ with $\omega(N)=3, 4, 5, 6, 7, 8$ and $9$ are 
$76546075$, $4784137$,  $20995$, $19699$, $2267$, $785$ and $299$, 
respectively.

\end{theorem}

To begin the proof of Theorem \ref{Vanishing}, observe that if equation (\ref{varnish}) has a vanishing subsum, then one of the following occurs :

\begin{eqnarray}
2^{a_1} 3^{b_1} + 2^{c_1} - 2^{a_2} 3^{b_2} = 0, \; \; 3^{d_1} - 2^{c_2} - 3^{d_2} = 0,  \label{eq1} \\
2^{a_1} 3^{b_1} + 2^{c_1} - 2^{c_2} = 0, \; \; 3^{d_1} - 2^{a_2} 3^{b_2} - 3^{d_2} = 0,  \label{eq2} \\
2^{a_1} 3^{b_1} + 2^{c_1} - 3^{d_2}=0, \; \; 3^{d_1} - 2^{a_2} 3^{b_2} - 2^{c_2} = 0,  \label{eq3} \\
2^{a_1} 3^{b_1} + 3^{d_1} - 2^{a_2} 3^{b_2} = 0, \; \; 2^{c_1} -2^{c_2} - 3^{d_2} = 0, \label{eq4} \\
2^{a_1} 3^{b_1} + 3^{d_1} - 2^{c_2} = 0, \; \; 2^{c_1} - 2^{a_2} 3^{b_2} - 3^{d_2} = 0,  \label{eq5} \\
2^{c_1}+3^{d_1} - 2^{a_2} 3^{b_2} = 0, \; \; 2^{a_1} 3^{b_1} - 2^{c_2} - 3^{d_2} = 0, \label{eq6} \\
2^{a_1} 3^{b_1} -2^{a_2} 3^{b_2}=0, \; \; 2^{c_1}+3^{d_1} - 2^{c_2} - 3^{d_2}=0,  \label{eq7} \\
2^{c_1}-2^{c_2}=0, \; \; 2^{a_1} 3^{b_1} +3^{d_1}-  2^{a_2} 3^{b_2} - 3^{d_2}=0,  \label{eq8} \\
3^{d_1}-3^{d_2}=0, \; \; 2^{a_1} 3^{b_1} +2^{c_1}-  2^{a_2} 3^{b_2} - 2^{c_2}=0. \label{eq9}
\end{eqnarray}

We consider each of these pairs of equations in turn; this argument is carried out in Tijdeman and Wang, essentially by reducing the problem to Propositions \ref{3-termite} and \ref{4-termite}. We will need to be somewhat more explicit and to appeal additionally to Theorem \ref{MT2}.

 In case (\ref{eq1}), necessarily $d_2=0$ and hence, from Propositions \ref{3-termite} and \ref{4-termite},  
$$
(d_1,c_2,d_2)= (1,1,0) \mbox{ or } (2,3,0),
$$
and $(a_1,b_1,c_1,a_2,b_2) $ is one of 
$$
\begin{array}{c} 
 (a_1,0,a_1,a_1+1,0), (a_1,1,a_1,a_1+2,0), (a_1,0,a_1-1,a_1-1,1), \\
 (a_1,0,a_1-3,a_1-3,2), (a_1,0,a_1+1,a_1,1), (a_1,0,a_1+3,a_1,2).
 \end{array}
$$
We thus have, in each case, that either $N$ is type I special or that (\ref{type3-3}) is satisfied, so that $N$ is type III special.

Arguing similarly, we find that 
in case (\ref{eq2}), $N$ is necessarily type I special. 
In cases (\ref{eq3}), (\ref{eq5}) and (\ref{eq6}), we have $N \leq 18$, while, in cases (\ref{eq4}), (\ref{eq7}), (\ref{eq8}) and (\ref{eq9}), $N$ is either type I special, or 
we have equations (\ref{type3-4}),  (\ref{type3-5}), either (\ref{type2-1}) or  (\ref{type3-1}), and either (\ref{type2-2}) or  (\ref{type3-2}), respectively.

If $N$ is type I special, or if equations (\ref{type3-3}) or (\ref{type3-4}) are satisfied, then we can appeal to  Theorem \ref{MT2} to solve the equation
$$
N = 2^x3^y+2^z+3^w
$$
without vanishing subsums; Theorem \ref{Vanishing} is immediate in these cases. To handle the remaining possibilities, we argue as in the proof of Theorem \ref{MT2}. By way of example,  consider the equation
\begin{equation} \label{cuty}
N=2^a+3^bc= 2^x3^y+2^z+3^w, \; \; c \in \{ 5, 7, 11, 13, 17, 19, 25, 35, 43, 73, 97, 145, 259 \},
\end{equation}
where we assume that
\begin{equation} \label{arggh}
2^a \not\in \{ 2^x3^y, 2^z, 3^w, 2^x3^y+2^z, 2^x3^y+3^w, 2^z+3^w \}. 
\end{equation} 
We have that
\begin{equation} \label{starterpack}
2^{\min \{ a, x, z \}} \mid 3^b c-3^w \; \; \mbox{ and } \; \; 3^{\min \{ b, y, w \} } \mid 2^a-2^z.
\end{equation}
Write
$$
M = \max \{ 2^a, 3^b c \} \; \; \mbox{ and } \; \; 
\{ 2^x 3^y, 2^z, 3^w \} = \{ M_1, M_2, M_3 \},
$$
as multisets, with $M_1 \geq M_2 \geq M_3$.  Suppose initially that $M$ satisfies 
\begin{equation} \label{M-lower-new}
\log M > 10^{20}.
\end{equation}
 
Let us assume first that $M=2^a$. Then we may apply  Proposition \ref{LFL2} to $M-M_1=M_2+M_3-  3^b c$ to conclude that
\begin{equation} \label{frenchies}
\max \{ M_2, 3^b c  \} > \frac{M}{3}  \cdot \mbox{exp} \left( -20 \log^2 (\log M) \right).
\end{equation}
 If, further, $M_1=2^z$, then $z=a-1$ and so
$$
2^{a-1} +3^bc= 2^x3^y+3^w,
$$
whereby $\min \{ b, y, w \} = 0$.
Writing in each case $c=2^\alpha+3^\beta$, for $\alpha$ and $\beta$ nonnegative integers, with
$$
(\alpha,\beta) \in \{ (2,0), (2,1), (1,2), (2,3), (4,0), (4,1), (4,2), (3,3), (4,3), (6,2), (4,4), (6,4), (4,5)  \},
$$
we thus obtain one of equations (\ref{eq07}), (\ref{eq14}) or (\ref{eq16}) and can once again appeal to Theorem \ref{MT2} to conclude that $N \leq 16481$.

We may thus suppose that $M_1 \in \{ 2^x3^y, 3^w \}$. If $\min \{ a, x, z \}=x$ and $\min \{ b, y, w \}=b$, then we match $2^x3^y$ and $3^w$, if $M_3=2^z,$ or $2^a$ and $2^z$ (if $M_2=2^z$), 
writing $2^x3^y+3^w=3^{\min \{ y, w \}} \beta$ or $2^a-2^z=2^z \beta$, respectively. In the first case, from  (\ref{starterpack}),   (\ref{M-lower-new}),  (\ref{frenchies}) and Proposition \ref{LFLP}, 
\begin{equation} \label{bet-bound}
0 < \beta < 2^{17 \left( \log \log M \right)^2 \log 259} \left( 1+ 3\cdot  e^{\left( \log \log M \right)^2} \right)< e^{67 \left( \log \log M \right)^2}.
\end{equation}
But then applying Theorem \ref{matveev} with $n=3$, we find  that 
$$
\log \left| \min \{ y, w \}  \log 3 + \log \beta - a \log 2 \right| > -2.2 \times 10^{11} \log (e a) \log \beta,
$$
and so, since
$$
\left| 2^a - 3^{\min \{ y, w \}} \beta \right| > 2^{a-1} \left| \min \{ y, w \}  \log 3 + \log \beta - a \log 2 \right|,
$$
it follows that
$$
2^z >  \frac{2^{a-1}}{\exp \left( 2.2 \times 10^{11} \log (e a) \log \beta \right)}.
$$
Thus
\begin{equation} \label{ZA}
z >  a-1- \frac{2.2 \times 10^{11} \log (e a) \log \beta}{\log 2} > 0.9 a.
\end{equation}
where the last inequality is a consequence of  (\ref{M-lower-new}) and  (\ref{bet-bound}). Applying Proposition \ref{LFLP} to 
$$
3^{\min \{ y, w \}-b} \beta -c,
$$ 
it follows that
$$
z < 17 \log^2 (\min \{ y, w \})  \log \beta,
$$
contradicting  (\ref{M-lower-new}), (\ref{bet-bound}) and (\ref{ZA}). If, on the other hand, $M_2=2^z$, so that $2^a-2^z=2^z \beta$ with
\begin{equation} \label{bet-too}
0 < \beta = 2^{a-z}-1 < \frac{\log M}{\log 2},
\end{equation}
we have 
$$
\left| 2^z \beta - \max \{ 2^x 3^y, 3^w \} \right| > \frac{M}{6\exp \left( 2.2 \times 10^{11} \log (e z) \log \beta \right)},
$$
whence
$$
\min \{ 2^x 3^y, 3^w \} > \frac{M}{6\exp \left( 2.2 \times 10^{11} \log (e z) \log \beta \right)}.
$$
Matching $2^x 3^y$ and $3^w$,  this contradicts (\ref{M-lower-new}), (\ref{bet-too}) and Proposition \ref{LFLP}.

Similarly, if $\min \{ a, x, z \}=x$ and $\min \{ b, y, w \}=y$ (so that $M_1=3^w$), we may apply  Proposition \ref{LFL2} to $2^a-3^w$ to conclude that
\begin{equation} \label{frenchies2}
 \max \{ 3^b c, 2^z \}   > M \cdot \mbox{exp} \left( -20 \log^2 (\log M) \right).
\end{equation}
If $M_3=2^z$, since $3^b \mid 2^a-\gamma$ where 
$$
0 < \gamma=2^x3^y+2^z \leq 2^{x+1} 3^y,
$$
this, with (\ref{starterpack}),  contradicts Proposition \ref{LFLP2}. If $M_2=2^z$ and $3^b c \geq 2^z$, then we match $3^bc$ and $3^w$, writing $3^w-3^bc=3^b \beta$. Using Theorem \ref{matveev} with $n=3$ to find a lower bound upon $|2^a-3^b \beta|$, this then implies a lower bound for $z$ that, with  (\ref{M-lower-new}), contradicts Proposition \ref{LFLP}.  If $M_2=2^z>3^b c$, 
then we match $2^a$ and $2^z$, apply  Theorem \ref{matveev} with $n=3$ to $|2^z(2^{a-z}-1) - 3^w|$ to derive a strong lower bound upon $3^bc$, match $3^w$ and $3^bc$, and arrive at a contradiction upon combining  (\ref{M-lower-new}) with Proposition \ref{LFLP2}.

Our arguments in the remaining cases (including those where $3^bc > 2^a$) are similar. After all is said and done, we are thus left to handle the values of $M$ with
 \begin{equation} \label{M-upper}
\log M \leq 10^{20}.
\end{equation}
As in the previous section, we consider cases separately. By way of example, if 
$$
\min \{ a, x, z \} = x, \; \; \min \{ b, y, w \} = w, \; \; M=3^bc, \; \; (M_1,M_2,M_3)=(2^x3^y,2^z,3^w),
$$
then (\ref{starterpack}) and (\ref{M-upper}) imply that $w \leq 42$. For each choice of $c$, we use Maple's ``ModularLog'' function to verify that the smallest positive integer $b-w$ for which
$$
3^{b-w} c \equiv 1 \mod{2^{74}}
$$
either does not exists or exceeds $10^{20}$, in each case, with the smallest value $105567241831233586666$ corresponding to $c=25$. It follows from (\ref{starterpack}) and (\ref{M-upper})  that $x \leq 73$.

Using
$$
3^{b-1} c \leq 2^x 3^y \leq 3^{b+1}c,
$$
we have 
$$
\frac{c}{3 \cdot 2^x} \leq 3^{y-b} \leq \frac{3c}{2^x},
$$
and hence for each pair $(x,c)$, sharp bounds upon $y-b$ and, in particular, leaves us with precisely $1924$ triples $(x,c,y-b)$ to consider. Writing
$$
3^bc-2^x3^y=3^{\min \{b,y \}} \beta,
$$
and assuming that $\min \{ b, y \} \geq w$, we thus have that $2^{\min \{ a, z \}} \mid 3^{\min \{ b, y \} -w} \beta-1$. For each of the $1792$ choices for $\beta$, we verify via Maple's ``ModularLog'' function that, in each case
$$
3^{\min \{ b, y \} -w} \beta \equiv 1 \mod{2^{81}}
$$
implies that $\min \{ b, y \} -w > 10^{20}$. We thus have $\min \{ a, z \} \leq 80$. We next bound $\max \{ a, z \}$ by considering 
$$
3^{\min \{b,y \}} \beta - 3^w \pm 2^{\min \{ a, z \}} \mod{2^{90}}.
$$
For each of our
$$
1792 \times 43 \times 81 \times 2 = 12483072
$$
cases, we find that necessarily $\max \{ a, z \} \leq 91$. With these upper bounds in hand, it then follows from 
$$
2^z-2^a+3^w =3^{\min \{b,y \}} \beta,
$$
that, crudely, $\min \{b, y \} \leq 56$.

To finish our analysis of this case, we simply find all solutions of equation (\ref{cuty}) -- with or without vanishing subsums -- modulo the integer 
\[
N_{180} = \gcd(2^{180}-1,3^{180}-1) = 439564261361225.
\]
This is computationally fairly straightforward, especially if we first compute solutions modulo 
\[
N_{36} = \gcd(2^{36}-1,3^{36}-1) = 23350145.
\]
Since our exponents satisfy $\max \{ a, b, x,y,z,w \} < 180$, it is a routine matter to check whether local solutions correspond to  actual solutions. In all cases, we find that solutions without vanishing subsums have
$$
a \leq 21, \; \; b \leq 10, \; \; x \leq 18, \; \; y \leq 10, \; \; z \leq 21 \; \mbox{ and } \; w \leq 13.
$$ 
The largest value of $N=2^a+3^bc$ with a representation as $2^x3^y+2^z+3^w$ with (\ref{arggh}) is $N=2099483$.
This completes the proof of Theorem \ref{Vanishing} for this case. The other cases proceed in a similar fashion.

\section{Newman's Conjecture : six-term equations  and beyond} \label{newman3}

In light of Theorem \ref{Vanishing}, to complete the proof of Theorem \ref{thm-New}, it suffices to consider positive integers $N$
for which $\omega(N) \geq 2$, say
$$
N=2^{a_1}3^{b_1} + 2^{c_1} + 3^{d_1} = 2^{a_2}3^{b_2} + 2^{c_2} + 3^{d_2},
$$
with the $a_i, b_i, c_i$ and $d_i$  nonnegative integers for which
\[
\{ 2^{a_1}3^{b_1}, 2^{c_1}, 3^{d_1}  \} \neq \{ 2^{a_2}3^{b_2}, 2^{c_2}, 3^{d_2}  \},
\]
where the corresponding $S$-unit equation
$$
2^{a_1}3^{b_1} + 2^{c_1} + 3^{d_1} - 2^{a_2}3^{b_2} - 2^{c_2} - 3^{d_2}=0
$$
has no vanishing subsums.
From work of  Evertse \cite{Ev} (see also van der Poorten and Schlickewei \cite{PS}), this equation has at most finitely many solutions. Unfortunately, we are unable to deduce an effective analogue of this statement. We will instead prove the following.

\begin{theorem} \label{non-Vanishing}
If $N$ is a positive integer with $\omega (N) \geq 3$ and $N$ is not special of type I, II or III,
then
\begin{equation} \label{jet}
N \in \{ 274, 473, 505, 1109, 1595, 1811, 2297, 2779, 4403, 20761 \}.
\end{equation}
For each of these values, we have $\omega (N)=3$.
\end{theorem}

Let us proceed by assuming, then, that $\omega (N) \geq 3$, so that
\begin{equation}\label{threereps}
N = 2^{a_1}3^{b_1} + 2^{c_1} + 3^{d_1} = 2^{a_2}3^{b_2} + 2^{c_2} + 3^{d_2} = 2^{a_3}3^{b_3} + 2^{c_3} + 3^{d_3},
\end{equation}
where again the $a_i, b_i, c_i$ and $d_i$ are nonnegative integers and
 \[
\{ 2^{a_i}3^{b_i}, 2^{c_i}, 3^{d_i}  \} \neq \{ 2^{a_j}3^{b_j}, 2^{c_j}, 3^{d_j}  \}
\]
for $i \neq j$, and, in each case, we will suppose that there are no vanishing subsums among the three resulting equations given by
\begin{equation} \label{fisher-king}
2^{a_i}3^{b_i} + 2^{c_i} + 3^{d_i} - 2^{a_j}3^{b_j} - 2^{c_j} - 3^{d_j}=0, \; \; 1 \leq i < j \leq 3.
\end{equation}
Our basic goal will be to reduce one of these three equations to a five-term equation to which we can apply Theorem \ref{eqmainthm1}. 
As before, note that, if $k$ is a positive integer, we have
\[
\nu_2(3^k-1) \le \nu_2(k) +2 \quad \text{and} \quad \nu_3(2^k-1)\le \nu_3(k)+1.
\]
Since we are assuming that there are no vanishing subsums, it follows from (\ref{threereps})  that
\begin{equation}\label{smallexp}
\min\{a_i,a_j,c_i,c_j\} \le \frac{\log|d_i-d_j|}{\log 2}+2 \; \; \mbox{ and } \; \;  \min\{b_i,b_j,d_i,d_j\}\le \frac{\log|c_i-c_j|}{\log 3}+1,
\end{equation}
for $1\le i,j\le 3, i\neq j$, whereby
\begin{equation}\label{smallterm}
2^{\min\{a_i,a_j,c_i,c_j\}} \le \frac{4\log N}{\log3}\; \; \mbox{ and } \; \;  3^{\min\{b_i,b_j,d_i,d_j\}}\le \frac{3\log N}{\log 2}.
\end{equation}
We will, by abuse of notation, refer to the quantities
$$
2^{a_i}, 3^{b_i}, 2^{c_i}, 3^{d_i}, \; \; 1 \leq i \leq 3
$$
as {\it monomials} and to 
$$
2^{a_i}3^{b_i}, 2^{c_i}, 3^{d_i}, \; \; 1 \leq i \leq 3
$$
as {\it terms}.
From (\ref{smallterm}), it follows that four of the monomials  $2^{\alpha_1}, 2^{\alpha_2}, 3^{\beta_1}, 3^{\beta_2}$ in the set
$$
\{ 2^{a_i}, 3^{b_i}, 2^{c_i}, 3^{d_i} \; : \; 1 \leq i \leq 3 \}
$$
are {\it small}  (in that their corresponding exponents are bounded by a constant multiple of $\log \log N$). Since the monomials $2^{\alpha_1}$ and $2^{\alpha_2}$ correspond to different representations of $N$ (as do $3^{\beta_1}$ and $3^{\beta_2}$), we may conclude  that we either have two representations of $N$ with two  small  monomials each, or one representation with two small monomials, while the other two representations have one small  monomial each.

We will show that it is always possible to match two terms from at least one pair of equations in \eqref{threereps}, and thus obtain  a $5$-term $S$-unit equation that we can solve effectively, via  Theorem  \ref{eqmainthm1}.
Note that a representation of $N=2^{a_i}3^{b_i} + 2^{c_i} + 3^{d_i}$ having two  small monomials immediately implies the existence of a small term with
$$
\min \{ 2^{a_i}3^{b_i},  2^{c_i},  3^{d_i} \} \leq \frac{4\log N}{\log3} \times \frac{3\log N}{\log 2} < 16 \log^2 N.
$$
If we have two representations of $N$, each with two small monomials (and hence at least one small term each), the corresponding equation (\ref{fisher-king}) can thus 
be rewritten as
\begin{equation} \label{easy-street}
u_1+u_2+u_3+u_4+\alpha=0,
\end{equation}
where the $u_i$ are $\{2, 3 \}$-units and $\alpha$ is an integer (the sum or difference of the two small $S$-units) with $|\alpha| < 16 \log^2N$. Such an equation may be effectively solved through appeal  to Theorem \ref{eqmainthm1}.

Suppose then that we have one representation of $N$, say   $N=2^{a_1}3^{b_1} + 2^{c_1} + 3^{d_1}$, with two small monomials, and two other representations with one  small monomial each. We may assume that
$$
\min \{ 2^{a_2}3^{b_2},  2^{c_2},  3^{d_2} , 2^{a_3}3^{b_3},  2^{c_3},  3^{d_3}  \} \geq 16 \log^2 N,
$$
or else reduce again to equation (\ref{easy-street}) with $|\alpha| < 16 \log^2N$. Without loss of generality, we have
$$
2^{a_2} \leq  \frac{4\log N}{\log3} \; \; \mbox{ and } \; \; 3^{b_3} \leq  \frac{3\log N}{\log2}.
$$
Note that
$$
\max \{ 2^{a_i} 3^{b_i}, 2^{c_i}, 3^{d_i} \} > N/3.
$$
Since we assume that there are no vanishing subsums,  the $c_i$ and $d_i$ are necessarily distinct so that, in particular, $3^{d_i}> N/3$ for at most one index $i$ and $2^{c_i} > N/3$ for at most two indices which must differ by $1$.
If we have either
$$
\max \{ 2^{a_2} 3^{b_2}, 2^{c_2}, 3^{d_2} \} =  2^{a_2} 3^{b_2} \; \; \mbox{ and } \; \; \max \{ 2^{a_3} 3^{b_3}, 2^{c_3}, 3^{d_3} \} = 3^{d_3}
$$
or
$$
\max \{ 2^{a_2} 3^{b_2}, 2^{c_2}, 3^{d_2} \} =2^{c_2}  \; \; \mbox{ and } \; \; \max \{ 2^{a_3} 3^{b_3}, 2^{c_3}, 3^{d_3} \} \in \{ 2^{a_3} 3^{b_3}, 2^{c_3} \},
$$
then  we can match terms in equation (\ref{fisher-king})  with $(i,j)=(2,3)$ to obtain an equation of the shape (\ref{easy-street}) with
$$
\alpha = 
\begin{cases}
3^{b_2} \beta, \; \; \beta = 2^{a_2} - 3^{d_3-b_2}, \; \; |\beta| < \frac{8\log N}{\log3}, \; \mbox{ if } 2^{a_2} 3^{b_2}, 3^{d_3} > N/3, \\
\min \{ 2^{c_2}, 2^{c_3} \},  \; \mbox{ if } 2^{c_2}, 2^{c_3} > N/3, \\
\min \{ 2^{c_2}, 2^{a_3} \},  \; \mbox{ if } 2^{c_2}, 2^{a_3} > N/3, \\
2^{a_3} \beta, \; \; \beta = 3^{b_3} - 2^{c_2-a_3}, \; \; |\beta| < \frac{3\log N}{\log2}, \; \mbox{ if } 2^{a_3} 3^{b_3}, 2^{c_2} > N/3, \; b_3 \geq 1; \\
\end{cases}
$$
in each case, effectively solvable via Theorem \ref{eqmainthm1}.
Otherwise, we necessarily have
either
$$
\max \{ 2^{a_2} 3^{b_2}, 2^{c_2}, 3^{d_2} \} \in \{   2^{a_2} 3^{b_2}, 3^{d_2} \}  \; \; \mbox{ and } \; \; \max \{ 2^{a_3} 3^{b_3}, 2^{c_3}, 3^{d_3} \} \in \{ 2^{a_3} 3^{b_3},  2^{c_3} \},
$$
or
$$
\max \{ 2^{a_2} 3^{b_2}, 2^{c_2}, 3^{d_2} \} =2^{c_2}  \; \; \mbox{ and } \; \; \max \{ 2^{a_3} 3^{b_3}, 2^{c_3}, 3^{d_3} \}=3^{d_3}.
$$
If 
$$
2^{c_1} \leq  \frac{4\log N}{\log3} \; \; \mbox{ and } \; \; 3^{d_1} \leq  \frac{3\log N}{\log2},
$$
then we reach an equation of the shape  (\ref{easy-street}) with $\alpha = 2^{c_1}+3^{d_1} < 8 \log N$.  Otherwise, since the representation  $N=2^{a_1}3^{b_1} + 2^{c_1} + 3^{d_1}$ has two small  monomials $2^e$ and $3^f$, we can match the term $\max \{ 2^{a_1} 3^{b_1}, 2^{c_1}, 3^{d_1} \}$ with either $\max \{ 2^{a_2} 3^{b_2}, 2^{c_2}, 3^{d_2} \} $ or $\max \{ 2^{a_3} 3^{b_3}, 2^{c_3}, 3^{d_3} \}$ in equation (\ref{fisher-king}) with either $(i,j)=(1,3)$ or $(i,j)=(2,3)$  to reduce to a five-term equation to which we can apply Theorem \ref{eqmainthm1}.

As is apparent, to carry this argument out in a completely explicit fashion, one is led to a rather involved case-by-case analysis. It is also by no means obvious that the bounds we achieve will be small enough to permit us to explicitly solve the resulting equations. We will provide full details for the case that is the most involved computationally, leading to the largest upper bounds upon our exponents,  and including all the features of our approach. Specifically, we deal with the case where the three equations in (\ref{fisher-king}) have no vanishing subsums and 
\begin{equation} \label{Caser1}
\begin{array}{c}
\min\{a_1,a_2,c_1,c_2\} = a_1, \; \min\{a_2,a_3,c_2,c_3\} = a_3, \;  \min\{b_1,b_3,d_1,d_3\} = d_1, \;  \\
\min\{b_2,b_3,d_2,d_3\} = b_2, \;  \max \{ 2^{a_1} 3^{b_1}, 2^{c_1}, 3^{d_1} \} = 2^{a_1}3^{b_1}, \;  \max \{ 2^{a_3} 3^{b_3}, 2^{c_3}, 3^{d_3} \} =  2^{a_3} 3^{b_3}. \\
\end{array}
\end{equation}

To treat  (\ref{Caser1}), we begin by supposing that
\begin{equation} \label{big-bound}
\log N  \geq 1.4 \times 10^{28}.
\end{equation}
We consider equation (\ref{fisher-king}) with $(i,j)=(1,3)$ and match the terms $2^{a_1}3^{b_1}$ and $2^{a_3} 3^{b_3}$, noting that we cannot have $a_1=a_3$. If $a_1 >a_3$, since 
$$
\min \{ 2^{a_1}3^{b_1}, 2^{a_3} 3^{b_3} \} > N/3,
$$
necessarily $b_3 \geq b_1$ and hence we can write
$$
2^{a_1}3^{b_1}-2^{a_3} 3^{b_3} = 2^{a_3} 3^{b_1} \left( 2^{a_1-a_3} - 3^{b_3-b_1} \right) = 2^{a_3} 3^{b_1} \beta,
$$
where $\beta \in \mathbb{Z}$. If $\beta > 0$, then 
$$
\beta < 2^{a_1-a_3} < \frac{4\log N}{\log3}.
$$
If, on the other hand, $\beta < 0$, since $2^{a_3} 3^{b_3} < 2^{a_1}3^{b_1+1}$, we have
\begin{equation} \label{better-bound}
|\beta| < 3 \cdot 2^{a_1-a_3} - 2^{a_1-a_3} = 2^{a_1-a_3+1} <\frac{8\log N}{\log3}.
\end{equation}
Similarly, if $a_3 > a_1$, $b_1 \geq b_3$ and we can write $2^{a_1}3^{b_1}-2^{a_3} 3^{b_3} =2^{a_1} 3^{b_3} \beta$, where again $|\beta| < \frac{8\log N}{\log3}$.
Equation (\ref{fisher-king}) with $(i,j)=(1,3)$ thus becomes
\begin{equation} \label{progress}
2^{c_1}+3^{d_1} + 2^a 3^b \beta - 2^{c_3}-3^{d_3}=0,
\end{equation}
where
\begin{equation} \label{conditions}
2^a \leq \frac{4\log N}{\log3}, \; \; 3^{d_1} \leq \frac{3 \log N}{\log2}, \; \; |\beta| < \frac{8\log N}{\log3}\; \; \mbox{ and }  \; \; 3^{b} > \frac{N\log 3}{12 \log N}.
\end{equation}
Suppose that $\beta > 0$; the case with $\beta < 0$ is very similar. Then if 
$$
\max \{ 2^{c_1},  2^a 3^b \beta \} = 2^{c_1} \; \; \mbox{ and } \; \; \max \{ 2^{c_3}, 3^{d_3} \} = 2^{c_3},
$$
the fact that $d_3 > d_1$ contradicts $c_1 \neq c_3$. If
$$
\max \{ 2^{c_1},  2^a 3^b \beta \} = 2^{c_1} \; \; \mbox{ and } \; \; \max \{ 2^{c_3}, 3^{d_3} \} = 3^{d_3},
$$
then writing
$$
3^{d_3} - 2^a 3^b \beta = 3^b \left( 3^{d_3-b}-2^a \right) = 3^b \beta_1,
$$
we have, from (\ref{conditions}), that $0 < \beta_1 < \frac{12\log N}{\log3}$, and that
$$
2^{c_1}+3^{d_1}  - 2^{c_3}-3^{b} \beta_1=0,
$$
whence $2^{c_3} \mid \beta_1 3^{b-d_1}-1$. Appealing to Proposition \ref{LFLP}, it follows that either $b-d_1 \leq 10^5$, or that
$$
c_3 < 17 \log^2 (b-d_1) \log (\beta_1). 
$$
In the first case, from (\ref{conditions}), we find that
$$
 \log \left( \frac{N \log 2 \log 3}{36 \log^2N} \right) < 10^5 \log 3,
$$
whence $\log N < 1.1 \times 10^5$, contradicting (\ref{big-bound}). 
In the second case,
$$
c_3 < 17 \log^2 \left( \frac{\log N}{\log 3} \right) \log \left( \frac{12\log N}{\log3} \right)
$$
and hence (\ref{progress}) can be rewritten as 
$$
2^{c_1}+\beta_2-3^{b} \beta_1=0,
$$
where $\beta_2=3^{d_1}-2^{c_3}$ satisfies
\begin{equation} \label{boo1}
\log |\beta_2| < \frac{17}{\log 2}  \log^2 \left( \frac{\log N}{\log 3} \right) \log \left( \frac{12\log N}{\log3} \right).
\end{equation}
In the other direction, we can apply Theorem \ref{matveev} with $n=3$ to deduce that 
$$
\log \left| b \log 3 + \log \beta_1 - c_1 \log 2 \right| > -2.2 \times 10^{11} \log (e c_1) \log \beta_1,
$$
whereby
\begin{equation} \label{boo2}
|\beta_2| = |2^{c_1}-3^b \beta_1| >  2^{c_1-1} \left| b \log 3 + \log \beta_1 - c_1 \log 2  \right| > 
\frac{2^{c_1-1}}{\exp \left( 2.2 \times 10^{11} \log (e c_1) \log \beta_1 \right)}.
\end{equation}
From (\ref{conditions}) and $\max \{ 2^{c_1},  2^a 3^b \beta \} = 2^{c_1}$, we have
$$
c_1 > \frac{\log \left( \frac{N\log 3}{12 \log N} \right)}{\log 2}
$$
and hence, combining (\ref{boo1}), (\ref{boo2}) and $0 < \beta_1 < \frac{12\log N}{\log3}$,  we find that
$\log N < 3 \times 10^{14}$,
contradicting (\ref{big-bound}).

We may thus suppose that $\max \{ 2^{c_1},  2^a 3^b \beta \} = 2^a3^b \beta$. If $\max \{ 2^{c_3}, 3^{d_3} \} = 2^{c_3}$, we again apply Theorem \ref{matveev} with $n=3$ to deduce a lower bound upon
$$
\left| 2^a3^b \beta -2^{c_3} \right| = 2^a \left| 3^b \beta - 2^{c_3-a} \right|.
$$
Specifically, from (\ref{conditions}), 
$$
\log \left| 2^a3^b \beta -2^{c_3} \right|  >c_3 \log 2   -2.2 \times 10^{11} \log (e (c_3-a) ) \log \left( \frac{8\log N}{\log3} \right)  - \log 2.
$$
Since $2^a3^b \beta -2^{c_3} = 3^{d_3}-2^{c_1}-3^{d_1}$, it follows from (\ref{conditions}) that $\max \{ 2^{c_1}, 3^{d_3} \}$ is close in size to both $2^a3^b \beta$ and $2^{c_3}$. In particular, we can match either $2^{c_1}$ and $2^{c_3}$, or $3^{d_3}$ and $2^a3^b \beta$. To be precise, if $\max \{ 2^{c_1}, 3^{d_3} \}=2^{c_1}$, then
$$
2^{c_3}-2^{c_1} = 2^{c_1} \left(2^{c_3-c_1}-1 \right) = 2^{c_1} \beta_3,
$$
where
$$
0 < \beta_3 < \frac{3 \cdot 2^{c_3}}{|2^a3^b \beta -2^{c_3}|} < 6 \exp \left( 2.2 \times 10^{11} \log (e (c_3-a) ) \log \left( \frac{8\log N}{\log3} \right) \right).
$$
We thus have from (\ref{big-bound}) that
\begin{equation} \label{beta3}
\log \beta_3 < 2.21 \times 10^{11} \left( \log ( \log (N ))\right)^2.
\end{equation}
From
$$
3^{d_1} + 2^a 3^b \beta - 2^{c_1} \beta_3-3^{d_3}=0,
$$
since $a, d_3 > d_1$, it follows that $3^{d_1} \| \beta_3$ and hence, applying Proposition \ref{LFLP2}  to $2^{c_1} (\beta_3/3^{d_1}) -1$, either $c_1 \leq 10^5$, or we have
$$
\min \{b, d_3 \} -d_1 < 10 \log^2 c_1 \log (\beta_3/3^{d_1}).
$$
From (\ref{conditions} and (\ref{beta3}), it follows that
\begin{equation} \label{dee3}
d_3 < 10 \log^2 c_1 \log \beta_3 < 2.22 \times 10^{12}   \left( \log ( \log (N ))\right)^4
\end{equation}
and hence the difference
$$
\left| 3^{b-d_1} \beta - 2^{c_1-a}  (\beta_3/3^{d_1}) \right| = 2^{-a} (3^{d_3-d_1}-1)
$$
is small. Appealing once again to Theorem \ref{matveev} with $n=3$, we find from
$$
3^{d_3} >  2^{a-1} 3^{b} \beta \left| (b-d_1) \log 3 -(c_1-a) \log 2 + \log (3^{d_1} \beta/\beta_3)  \right|
$$
and (\ref{beta3}), that
$$
d_3 > \frac{1}{\log 3} \left( (a-1) \log 2 + b \log 3 + \log \beta  - 4.87 \times 10^{22}  \log (e(c_1-a))   \left( \log ( \log (N ))\right)^2  \right).
$$
Combining this with (\ref{conditions}), (\ref{dee3}) and $2^{c_1} < N$ contradicts (\ref{big-bound}). 

If, on the other hand, $\max \{ 2^{c_1}, 3^{d_3} \}=3^{d_3}$,  we write 
$$
2^a 3^b \beta - 3^{d_3} = 3^{\min \{ b, d_3 \}} \beta_4,
$$
and arguing as previously again, after a little work, contradict (\ref{big-bound}).

\subsection{Handling smaller values of $N$}

To treat the remaining smaller values of $N$, we must work somewhat harder than in the preceding sections, employing rather more sophisticated arguments. We may suppose that
\begin{equation} \label{init-bound3}
\log N < 1.4 \times 10^{28}.
\end{equation}
To begin, note that, from (\ref{smallterm}), {\ref{Caser1}) and (\ref{init-bound3}), we have 
$$
a_1 \leq 95, \; a_3 \leq 95, \; d_1 \leq 60 \; \mbox{ and }  \; b_2 \leq 60,
$$
and 
$$
 \max \{ 2^{a_1} 3^{b_1}, 2^{c_1}, 3^{d_1} \} = 2^{a_1}3^{b_1} 
\mbox{ and } \max \{ 2^{a_3} 3^{b_3}, 2^{c_3}, 3^{d_3} \} =  2^{a_3} 3^{b_3}. 
$$
It follows that
$$
2^{a_1}3^{b_1-1} < 2^{a_3} 3^{b_3} < 2^{a_1}3^{b_1+1}
$$
and hence, since we have no vanishing subsums, $a_1 \neq a_3$. We therefore have
$$
 \frac{\log 2}{\log 3} (a_1-a_3) -1  <  b_3-b_1 < \frac{\log 2}{\log 3} (a_1-a_3) +1
$$
and hence, for a given pair $(a_1,a_3)$, precisely two choices for $b_3-b_1$.

To proceed, we will appeal to a result of Baker and Davenport, as it appears in Lemma 5 of Dujella and Peth\H{o} \cite{DuPe} :
\begin{prop} \label{BaDa}
Suppose that $M$ is a positive integer, and that $\kappa$ and $\mu$ are irrational real numbers . Let $p/q$ be a convergent in the infinite simple continued fraction expansion of  $\kappa$, satisfying $q > 6M$, and let
$$
\epsilon = \| \mu q \| - M \cdot \| \kappa q \|,
$$
where $\| \cdot \|$ denotes the distance to the nearest integer. If $\epsilon > 0$, and $A, B$ are positive real numbers with $B>1$,  then there is no solution to the inequality 
$$
0 < m \kappa - n +\mu < A B^{-m},
$$
in integers $m$ and $n$ with
$$
\frac{\log (Aq/\epsilon)}{\log B} \leq m \leq M.
$$
\end{prop}

For our purposes, if, say, $a_1<a_3$ and $3^{b_1-b_3}>2^{a_3-a_1}$, from
$$
(3^{b_1-b_3}- 2^{a_3-a_1}) 3^{b_3}  - 2^{c_3-a_1} = \frac{3^{d_3}-2^{c_1}-3^{d_1}}{2^{a_1}},
$$
it follows that 
\begin{equation} \label{goldfish}
\frac{3^{d_3}-2^{c_1}-3^{d_1}}{2^{c_3}} > b_3 \log 3 + \log (3^{b_1-b_3}- 2^{a_3-a_1}) - (c_3-a_1)  \log 2.
\end{equation}
We note for future reference that if we set
$$
p=599362460113865624518687902158 \; \mbox{ and } \; q = 378155609259623725703103696043,
$$
then $p/q$ is a convergent in the infinite simple continued fraction expansion of $\log (3)/\log (2)$.
If $3^{b_1-b_3}- 2^{a_3-a_1} = 1$, since $q > 1.4 \times 10^{28} > b_3$, it follows from the theory of continued fractions that 
$$
|b_3 \log 3 - 2^{a_3-a_1}) - (c_3-a_1)  \log 2| > |q \log 3 - p \log 2| > \left( 3.18 \times 10^{30} \right)^{-1}
$$
If, on the other hand,  $3^{b_1-b_3}- 2^{a_3-a_1} > 1$, we apply Proposition \ref{BaDa} with $\kappa = \frac{\log 3}{\log 2}$, $\mu = \frac{ \log (3^{b_1-b_3}- 2^{a_3-a_1})}{\log 2}$,  
$M=2.02 \times 10^{28} >  \frac{\log N}{\log 2}$ and $q$ as above. We verify that $M \| q \tau \| < 0.0107$ and that, in each case $\epsilon > 0$. We find further that $q/\epsilon < 2.77 \times 10^{31}$, with the largest value corresponding to 
$a_3-a_1=68$ and $b_1-b_3=43$. It follows from Proposition \ref{BaDa} that 
$$
b_3 \log 3 + \log (3^{b_1-b_3}- 2^{a_3-a_1}) - (c_3-a_1)  \log 2 > \left( 2.77 \times 10^{31} \right)^{-1}
$$
and hence, in all cases, 
$$
2^{c_3} <  \left( 2.77 \times 10^{31} \right) 3^{d_3}.
$$
From
$$
2^{c_3}+3^{d_3}-2^{c_1}-3^{d_1} = (3^{b_1-b_3}- 2^{a_3-a_1}) 2^{a_1} 3^{b_3},
$$
we therefore have
\begin{equation} \label{crunch}
 \left( 2.77 \times 10^{31} +1 \right)3^{d_3} > (3^{b_1-b_3}- 2^{a_3-a_1}) 2^{a_1} 3^{b_3}.
\end{equation}
This provides a strong lower bound upon $d_3-b_3$. In the worst case (when $a_1=0$ and $a_3 \in \{ 1, 3 \}$), for example, we have $d_3 - b_3 \geq -65$. We obtain similar bounds in case $a_1>a_3$ or $3^{b_1-b_3}<2^{a_3-a_1}$. To finish, we simultaneously match the terms $2^{a_1} 3^{b_1}$, $2^{a_3} 3^{b_3}$ and $3^{d_3}$, write
$$
2^{a_1} 3^{b_1}-2^{a_3} 3^{b_3}-3^{d_3} = 3^{\min \{ b_1, b_3, d_3 \}} \beta.
$$
For example, if 
\begin{equation} \label{goober}
a_1<a_3, \; d_3 < b_3 \; \mbox{ and } \; 3^{b_1-b_3}>2^{a_3-a_1},
\end{equation}
then 
$$
\beta = (3^{b_1-b_3}- 2^{a_3-a_1}) 2^{a_1} 3^{b_3-d_3} -1
$$
which, with (\ref{crunch}), provides precisely $64134$ values for $\beta$, assuming (\ref{goober}).
Note that
\begin{equation} \label{gusto}
3^{\min \{ b_1, b_3, d_3 \}} \beta + 3^{d_1} = 2^{c_3}-2^{c_1}.
\end{equation}
Appealing to Maple's ``ModularLog'' function, we find that the congruences 
$$
3^x \beta  \equiv -1 \mod{2^{112}}
$$
have, for each possible choice of $\beta$, no positive  solutions $x$ with $x <  1.4 \times 10^{28}$. We may thus conclude that $\min \{ c_1, c_3 \} \leq 111$.
Finally, we consider the congruences
$$
2^x =2^{\min \{ c_1, c_3 \}} \pm 3^{d_1} \mod{3^{71}},
$$
finding that we necessarily have either $x \leq 2$ or $x > 2.02 \times 10^{28}$, whereby we may conclude that $\max \{ c_1, c_3 \} \leq 2$, or that $\min \{ b_1, b_3, d_3 \} \leq 70$.
In either case, we conclude, after a little work, that
\begin{equation} \label{finished}
\max \{ a_1, b_1, c_1, d_1 \} <180.
\end{equation}
As previously, we consider the general equation
$$
2^{a_1}3^{b_1} + 2^{c_1} + 3^{d_1} = 2^{a_2}3^{b_2} + 2^{c_2} + 3^{d_2} 
$$
modulo $N_{180}$, finding all solutions (once again, with and without vanishing subsums). We check to see which if any of these correspond to actual values of $N$ with, additionally a third representation, no corresponding  vanishing subsums, and satisfying (\ref{Caser1}). It turns out that there are no such values of $N$, in this case. carrying out this argument for the remaining cases completes the proof of Theorem \ref{non-Vanishing} (and hence of Theorem \ref{thm-New}).

\section{Concluding remarks} \label{newman5}

The results of Tijdeman and Wang (specifically Theorem 4 of \cite{TW})  actually deal with the more general problem of representing rational numbers $N$ as in (\ref{newmanrep}), where now the variables $a, b, c, d$, while still integers, are allowed to be negative. The arguments of the paper at hand apply also in this more general situation, though some of the technical details are more involved.


\end{document}